\documentclass[a4paper,11pt]{amsart}
\usepackage{longtable, stackrel, multicol, fancybox, tcolorbox}
\usepackage{bm,amsfonts,amsmath,amssymb,mathrsfs,bbm,amsthm}
\usepackage{graphicx, color,eurosym, eucal,palatino}
\usepackage[normalem]{ulem}

\usepackage{marginnote}
\usepackage{xcolor}

\makeatletter
\newcommand\footnoteref[1]{\protected@xdef\@thefnmark{\ref{#1}}\@footnotemark}
\makeatother

\newtheorem{theorem}{Theorem}[section]

\title[Topologically induced spectral behavior]{\textsc{Topologically induced spectral behavior: \\ the example of quantum graphs}}

\author[P. Exner]{Pavel Exner}
\address{Doppler Institute for Mathematical Physics and Applied Mathematics\\
Czech Technical University in Prague\\ B\v{r}ehov\'{a} 7, 11519 Prague, Czech Republic,
{\rm and}
Department of Theoretical Physics\\
Nuclear Physics Institute, Czech Academy of Sciences,
25068 \v{R}e\v{z}, Czech Republic\\
E-mail: {exner@ujf.cas.cz}
}

\begin{document}
\begin{abstract}
This review paper summarizes the contents of the talk given by the author at the 8th International Congress of Chinese Mathematicians. Using examples of Schr\"odinger operators on metric graphs, it is shown that a nontrivial topology of the configuration space can give rise to a rich variety of spectral types. In particular, it is shown that the spectrum may be of a pure point type or to have a Cantor structure. We also address the question about the number of open spectral gaps and show that it could be nonzero and finite. Finally, inspired by a recent attempt to model the anomalous Hall effect we analyze a vertex coupling which exhibits high-energy behavior determined by the vertex degree parity.
\end{abstract}

\maketitle

%%%%%%%%%%%%%%%%%%%%%%%%%%%%%%%%%%%%%%
\section{Introduction} \label{s:intro}

Nice thing about mathematics are numerous relations between its different parts manifesting its deep unity. This is true, in particular, about spectral properties of operators and the topology of the underlying space; I believe that there is no need to convince the reader that the latter has an impact on the former. Generally speaking, one expects that a nontrivial topology can give rise to a broader family of spectral types. My aim here is to illustrate this claim using the example of \emph{quantum graphs} which I will introduce in the following section. Put like that, of course, the claim is rather vague; to be specific, I am going to consider \emph{elliptic second-order} operators with \emph{periodic coefficients} commonly used in physics to describe crystals and other periodically structured materials.

Let us first recall how things look like in the standard Schr\"odinger operator theory. A typical example of the indicated class of operators is
 % -------------- %
\begin{equation} \label{perSO}
H = (-i\nabla-A(x))^* g(x) (-i\nabla-A(x)) + V(x)
\end{equation}
 % -------------- %
on $L^2(\mathbb{R}^d)$, where $d\in\mathbb{N}$ is the space dimension, $g$ is a positive $d\times d$ matrix-valued function and $A$ is a vector-valued magnetic potential. If the coefficients $g,A,V$ are periodic the spectrum of the operator \eqref{perSO} is found using \emph{Floquet method}: one proves that $H$ is unitarily equivalent to
 % -------------- %
$$
H' = \int_{Q^*} H(\theta)\, \mathrm{d}\theta\,.
$$
 % -------------- %
The fiber operators $H(\theta)$ here act on $L^2(Q)$, where $Q\subset\mathbb{R}^d$ is the period cell and $Q^*$ is the dual cell, or \emph{Brillouin zone} in the physicist's terminology. Examining the spectrum of $H(\theta)$ as a function of the quasimomentum $\theta$, one can prove, in particular, that the spectrum of $H$

\vspace{-.5em}

 % -------------- %
\begin{itemize}
 \setlength{\itemsep}{0pt}
 \item is absolutely continuous, and
 \item has a band-and-gap structure.
\end{itemize}
 % -------------- %

\vspace{-.5em}

\noindent The argument is based on the analyticity of the eigenvalues with respect to $\theta$. Its idea belongs to L.~Thomas \cite{Th73} who used in the case $A=0$ and \mbox{$g=I$}, for a review and a general result under weak regularity assumptions see \cite{BS00}. Moreover, the gap structure is influenced by the dimension of the configuration space:

\vspace{-.5em}

\begin{itemize}
 \setlength{\itemsep}{0pt}

 \item in the one-dimensional case the number of open gaps is infinite except for a particular class of potentials,

 \item on the contrary, in higher dimensions the \emph{Bethe-Sommerfeld conjecture}, nowadays verified for a wide class of interactions \cite{Pa08}, says that the number of open gaps is \emph{finite}.

\end{itemize}

\vspace{-.5em}

\noindent I am going to show that if the system in question is a quantum graph, any of the above properties may appear to be invalid.

%%%%%%%%%%%%%%%%%%%%%%%%%%%%%%%%%%%%%
\section{Quantum graphs} \label{s:qg}
\setcounter{equation}{0}

Let us now introduce the main object of our interest. We consider a \emph{metric graph} which is a collection of vertices and edges connected in accordance with a prescribed adjacency matrix. In contrast to the most common concept of a discrete graph, we suppose in addition that each edge is homothetic to a (finite or semi-infinite) interval (cf.~Fig.~1).
 % ------------- %
 \begin{figure} %\label{f:metric}
  \centering
    \includegraphics[trim=5cm 24.3cm 7cm 2cm, width=12cm]{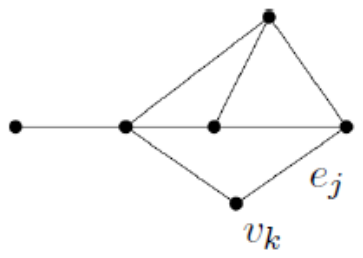}
 \setlength\unitlength{1mm}
 \caption{A metric graph}
 \end{figure}
 % ------------- %
We associate with such a metric graph the Hilbert space $\mathcal{H} = \bigoplus_j L^2(e_j)$ and consider the operator $H$ acting on functions $\psi=\{\psi_j\}$ that are locally $H^2$ as
 % -------------- %
$$
H\psi = \{-\psi''\} \quad\text{or more generally}\quad H\psi = \{ (-i\psi'-A\psi)^2 + V\psi\}
$$
 % -------------- %
Should such an operator play the role of a quantum mechanical Hamiltonian, it has to be self-adjoint. To achieve that we have to match the functions $\psi_j$ properly at each graph vertex.

Denoting $\psi(v_k)=\{\psi_j(v_k)\}$ and $\psi'(v_k)=\{\psi'_j(v_k)\}$ the  boundary values of functions and (outward) derivatives at the vertex $v_k$, respectively, the most general self-adjoint matching conditions read
 % ------------- %
\begin{equation} \label{bc}
 (U-I)\psi(v_k)+i(U+I)\psi'(v_k)=0,
\end{equation}
 % ------------- %
where $U$ is \emph{any} $\mathrm{deg}(v_k) \times \mathrm{deg}(v_k)$ unitary matrix. This is easy to see: an elementary calculation gives
 %------------%
 $$
 \|\psi(v_k)+i\psi'(v_k)\|^2 - \|\psi(v_k)-i\psi'(v_k)\|^2 =
 2\sum_j \big(\bar\psi_j \psi'_j - \bar\psi'_j \psi_j\big)(v_k)
 $$
 %-------------%
where the right-hand side is nothing but the boundary form that has to vanish to make the operator self-adjoint. Hence the above two vectors, $\psi(v_k)\pm i\psi'(v_k)\in \mathbb{C}^{\mathrm{deg}(v_k)}$, are of the same length being thus related by a unitary matrix.

In general the vertex coupling depends on the number of parameters, specifically $n^2$ of them for a vertex of degree $n$. The number is reduced significantly if we require \emph{continuity at the vertex}, then we are left with
 % ------------- %
\begin{equation} \label{delta}
 \psi_j(0)=\psi_k(0)=:\psi(0)\,,\;\; j,k=1,\dots,n\,,
 \;\;\;\;
 \sum_{j=1}^n \psi'_j(0) = \alpha \psi(0)
\end{equation}
 % ------------- %
depending on a single parameter $\alpha\in\mathbb{R}$. The conditions \eqref{delta} are called the $\delta$ \emph{coupling} and the corresponding unitary matrix is $U= {2\over n+i\alpha}\mathcal{J}-I$, where $\mathcal{J}$ is the $n\times n$ matrix whose all entries are equal to one. In particular, the case $\alpha=0$ is often referred to as the \emph{Kirchhoff coupling}. We note that this is an unfortunate name -- \emph{free} or \emph{standard} or \emph{natural} would be better -- but it stuck. The name $\delta$ coupling makes sense because one can approximate it by scaled regular potentials similarly as a $\delta$ potential on the line \cite{Ex96a}; the parameter $\alpha$ is at that interpreted as the coupling strength.

Before proceeding further, we have to say why quantum graphs are interesting. There are multiple reasons:

\vspace{-.5em}

\begin{itemize}

\item from the `practical' physical point of view they model graphlike \emph{nanostructures} made of semiconductors, carbon nanotubes, and other materials,

\item they offer a number of interesting mathematical questions, coming in the first place from the spectral and scattering theory, both direct and inverse (e.g., can on hear the shape of a graph? \cite{GS01}),

\item they have \emph{combinatorial graph} counterparts as we shall mention below, which may exhibit interesting dynamical effects, see, for instance \cite{LLY12},

\item graphs offer a useful playground to study \emph{quantum chaos} \cite{KS97},

\item and a lot more, as one can find, e.g., in the monograph \cite{BK13}.

\end{itemize}

We focus here, however, on a single aspect, namely how the topology can enrich spectral properties of quantum graphs.

%%%%%%%%%%%%%%%%%%%%%%%%%%%%%%%%%%%%%%%%%%%%%%%%%%%%%%
\section{The lack of absolute continuity} \label{s:ac}
\setcounter{equation}{0}

To describe how quantum graphs differ from the standard PDEs mentioned in the opening, we note first that the \emph{unique continuation principle} in general does not hold in quantum graphs, in other words, one can have compactly supported eigenfunctions. This is straightforward to see: a graph with a $\delta$ coupling at the vertices which contains a loop with commensurate edges has the so-called \emph{Dirichlet eigenvalues} corresponding to eigenfunctions of the type shown in Fig.~2 because the rest of the graph `does not know' about them.
 % -------------- %
\begin{figure}[h!]
\centering
    \includegraphics[scale=0.3,angle=0]{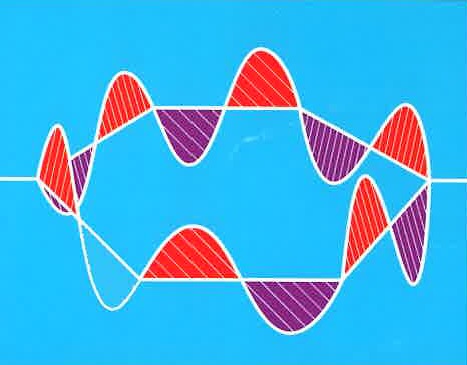}
    \caption{A Dirichlet eigenfunction; courtesy: Peter Kuchment}
\end{figure}
 % -------------- %
This means that the spectrum of a \emph{periodic} quantum graph with the said property contains infinitely degenerate eigenvalues being thus not purely absolutely continuous.

Other peculiarities of quantum graphs spectra are less trivial, nevertheless, they can be demonstrated using relatively simple examples. Next we want to show that such a spectrum may not have an absolutely continuous component at all; the example we use to this aim is a graph in the form of a \emph{loop array} exposed to a \emph{magnetic field} as sketched in Fig.~3; we suppose that the loops are circles of unit radii.
 % -------------- %
\begin{figure}[h!]
\centering
    \includegraphics[scale=0.8]{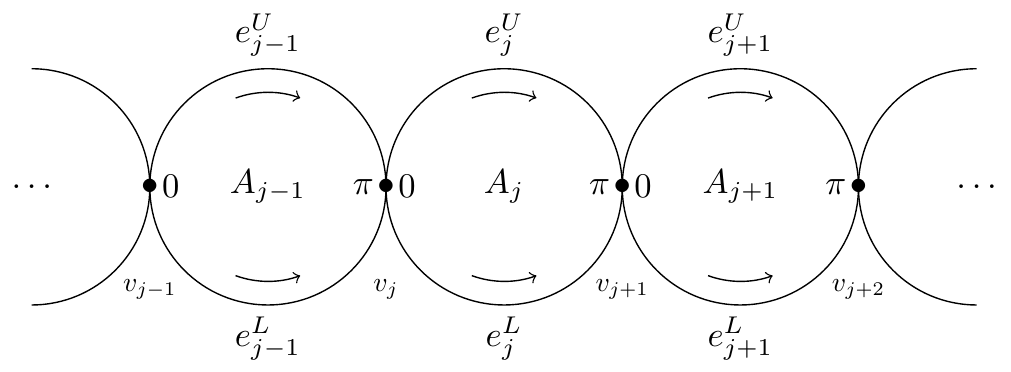}
     \caption{The magnetic chain graph; the arrows indicate the coordinate choice at the edges} %\label{fig:graph}
\end{figure}
 % -------------- %
The corresponding Hamiltonian is the magnetic Laplacian $-\Delta_A$ on the chain $\Gamma$ acting as $\psi_j\mapsto-\mathcal{D}^2\psi_j$ on each link of it, where $\mathcal{D}:= -i\nabla- \mathbf{A}$. We have some freedom in the choice of the vector potential. First of all, we need in fact its tangent component only, and secondly, it is only the flux through the loops which matters; without loss of generality we may suppose that the potential takes a constant value $A_j$ on the $j$-th loop. For the sake of definiteness we assume the magnetic version of the $\delta$-coupling at the vertices \cite{KS03}, in other words, the operator domain consists of functions from $H^2_\mathrm{loc}(\Gamma)$ satisfying
 % -------------- %
\begin{equation} \label{mg delta}
\psi_i(0) = \psi_j(0) = :\psi(0)\,,
\quad
i,j=1,\dots,n\,,
\quad
\sum_{j=1}^n \mathcal{D}\psi_j(0) =-i\alpha\,\psi(0)\,,
\end{equation}
 % -------------- %
where $\alpha\in\mathbb{R}$ is the coupling constant and $n=4$ holds in our case. In fact, the vector potential terms in this particular case cancel mutually at the left-hand side of the last of the conditions \eqref{mg delta} which thus looks as in \eqref{delta}). This does mean, however, that the magnetic potential does not affect the spectrum because it is also contained in the operator symbol -- and we will see in a while the magnetic effects may be quite dramatic.

Consider first the simplest situation with the homogeneous magnetic field. We use the circular gauge so that the (tangent component of) the vector potential is constant on each ring, $A_j=A,\: j\in\mathbb{Z}$, pointing clockwise. To find the spectrum we take the elementary cell of the chain as sketched in Fig.~4
 % -------------- %
\begin{figure}[h!]
\centering
    \includegraphics[scale=.8]{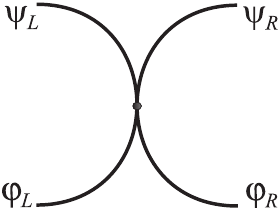}
     \caption{Elementary cell of the magnetic chain} %\label{fig:cell}
\end{figure}
 % -------------- %
and use the Ansatz $\psi_L(x)=\mathrm{e}^{-iAx} (C_L^+\mathrm{e}^{ikx}+C_L^-\mathrm{e}^{-ikx})$ for $x\in[-\pi/2,0]$ and energy $E:= k^2\neq0$, and similarly for the other three components; for $E<0$ we put instead $k=i\kappa$ with \mbox{$\kappa>0$}. These functions have to be matched through (a) the $\delta$-coupling and (b) the Floquet conditions at the `loose' ends. This yields a system of linear equations for the coefficients, its solvability condition is expressed as the following equation for the phase factor $\mathrm{e}^{i\vartheta}$,
 % -------------- %
$$
\sin k\pi\cos A\pi (\mathrm{e}^{2i\vartheta}-2\xi(k)\mathrm{e}^{i\vartheta}+1)=0\,,
$$
 % -------------- %
with the discriminant $D=4(\xi(k)^2-1)$, where $\xi(k) := \frac{\eta(k)}{4\cos A\pi}$ and
 % -------------- %
\begin{equation} \label{eta}
\eta(k) := 4\cos k\pi + \frac{\alpha}{k}\sin k\pi
\end{equation}
 % -------------- %
for any $k \in\mathbb{R}\cup i\mathbb{R}\setminus\{0\}$ and $A-\frac12\not\in\mathbb{Z}$. Putting aside the case with  $A-\frac12\in\mathbb{Z}$ and $k \in\mathbb{N}$, we conclude that $k^2\in\sigma(-\Delta_A)$ holds if and only if the condition $|\eta(k)| \leq 4|\cos A\pi|$ is satisfied. Its solution can be demonstrated graphically as shown in Fig.~5.
 % -------------- %
\begin{figure}
    \begin{center}
        \includegraphics[scale=0.9]{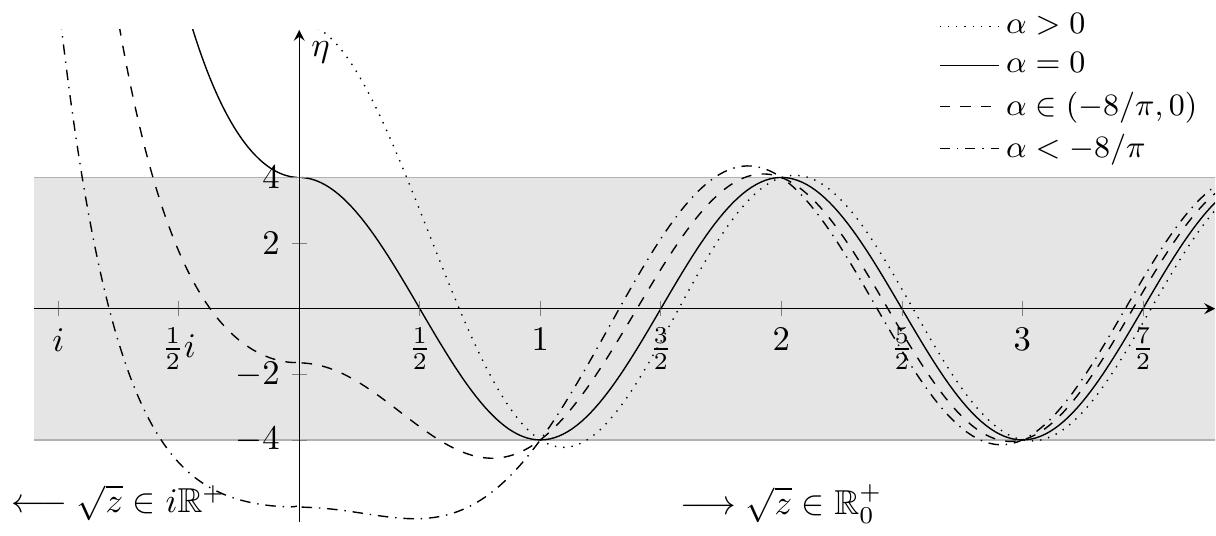} %% in pdfLaTeX
        \caption{Graphical solution of the magnetic chain spectral condition: the behavior of the quantity $\eta$ defined by \eqref{eta} with respect to the spectral parameter $z$ and the coupling constant $\alpha$ is shown. On the right side the horizontal axis indicates the variable $\sqrt{z}=k$ while on the left we plot increasing values of the purely imaginary positive values of the square root, $\sqrt{z} = i\kappa$,\, $\kappa > 0$.}
    \end{center}
    %\label{fig:eta_fun}
\end{figure}
 % -------------- %
The spectrum -- in the momentum variable -- consists of the values for which the values of the function \eqref{eta} lie in shaded strip. If $A=0$ and the coupling is Kirchhoff, the spectrum covers the positive real axis. A nonzero value of $\alpha$ opens spectral gaps, $\alpha<0$ means that the spectrum extends to negative values. If the magnetic field is nonzero -- or, in fact, non-integer -- the picture is similar but the strip width changes to $8|\cos A\pi|$ which means, in particular, that spectral gaps are open for any $\alpha\in\mathbb{R}$. What is important from our present point of view, for half-integer values, \mbox{$A-\frac12\in\mathbb{Z}$}, the strip shrinks to a line, and as a result, the spectrum consists of \emph{infinitely degenerate eigenvalues} only (or \emph{flat bands} as physicists would say), with the `elementary' eigenfunctions being supported by pairs of adjacent loops.

%%%%%%%%%%%%%%%%%%%%%%%%%%%%%%%%%%%%%%%%%%%%%%%%%
\section{Variable magnetic field} \label{s:varmg}
\setcounter{equation}{0}

Let us keep exploring the last example making it a little more complicated by assuming that the magnetic field is nonconstant and varies \emph{linearly} along the chain, $A_j = \mu j + \theta$ for some $\mu, \theta \in \mathbb{R}$ and every $j\in\mathbb{Z}$. Such a field extending over the whole space is an idealization, of course, but the validity of such idealizations is the bridge at which mathematics and physics meet at least since Newton times.

Being more pragmatic, we note that the unbounded character of the sequence $\{A_j\}$ need not bother us as it is not essential. The point is that from the spectral point of view only the \emph{fractional part of each $A_j$ matters.} The reason is that our operator -- which we denote as $-\Delta_{\alpha,A}$ for a given $\alpha\in\mathbb{R}$ and $A=\{A_j\} \subset \mathbb{R}$ -- is unitarily equivalent to $-\Delta_{\alpha,A'}$ with $A_j' = A_j + n$ with $n\in\mathbb{Z}$ by the operator acting as $\psi_j(x) \mapsto \psi_j(x)\, \mathrm{e}^{-inx}$; a physicist would call it a gauge transformation. This simplifies the analysis in the case when the slope $\mu$ is \emph{rational}. Indeed, is such a situation we can assume without loss of generality that the sequence $\{A_j\}$ is periodic and solve the problem using the Floquet method similarly as we did above for a constant $A$. In this way we arrive at the following result \cite{EV17}:

 % -------------- %
\begin{theorem}
  Let $A_j = \mu j + \theta$ for some $\mu, \theta \in \mathbb{R}$ and every $j\in\mathbb{Z}$. Then for the spectrum $\sigma(-\Delta_{\alpha,A})$ the following holds: \\[.3em]
 % -------------- %
  (a) If $\mu, \theta \in \mathbb{Z}$ and $\alpha = 0$, then the spectrum is the union of the absolutely continuous and pure point parts,
  $\sigma_{ac}(-\Delta_{\alpha,A}) = [0,\infty)$ and $\sigma_{pp}(-\Delta_{\alpha,A}) = \{n^2 |\, n\in\mathbb{N}\}$, \\[.15em]
 % -------------- %
  (b) If $\alpha \ne 0$ and $\mu = p/q$ with $p,q$ relatively prime, $\mu j + \theta + \frac{1}{2} \notin \mathbb{Z}$ for all $j = 0,\ldots, q-1$, then $-\Delta_{\alpha,A}$ has a family of infinitely degenerate eigenvalues, $\{n^2 |\, n\in\mathbb{N}\}$, interlaced with an absolutely continuous part consisting of $q$-tuples of closed intervals, \\[.15em]
 % -------------- %
  (c) If the situation is as in (b) but $\mu j + \theta + \frac{1}{2} \in \mathbb{Z}$ holds for some $j = 0,\ldots, q-1$, then the spectrum $\sigma(-\Delta_{\alpha,A})$ consists of infinitely degenerate eigenvalues only, the Dirichlet ones plus $q$ distinct others in each interval $(-\infty, 1)$ and $\big(n^2, (n+1)^2\big)$.
\end{theorem}

The case of an \emph{irrational $\mu$} requires a different approach: we shall rephrase our differential operator problem on the metric graph in terms of a \emph{difference equation}. This idea was put forward by physicists, Alexander and de Gennes, in the 1980s; its mathematical treatment started a decade later. It is particularly simple if the graph in question is \emph{equilateral} as is the case in our example. We introduce the set $\mathfrak{K}:=\{k\colon\mathrm{Im}\,k\geq0\wedge k\notin\mathbb{Z}\}$ with the aim to exclude Dirichlet eigenvalues from the consideration and seek the spectrum through solution of the Schr\"odinger equation {\scriptsize $(-\Delta_{\alpha,A}- k^2)\Big(\!\! \begin{array}{c} \psi(x,k)\\[-.2em] \varphi(x,k) \end{array} \!\!\Big) = 0$} in which the component functions refer to the upper and lower branch of the chain, respectively. They satisfy the condition \eqref{delta} at the vertices, $x=j\pi$, which allows us to reduce the problem to the difference equation
 % -------------- %
$$
2\cos(A_j\pi)\psi_{j+1}(k)
+
2\cos(A_{j-1}\pi)\psi_{j-1}(k)
=
\eta(k)\psi_j(k)\,,\;\;  k\in\mathfrak{K}\,,
$$
 % -------------- %
where $\psi_j(k):=\psi(j\pi,k)$ and $\eta$ is the function \eqref{eta} amended by $\eta(k)=4+\alpha\pi$ for $k=0$; note that since the upper and lower function component coincide at the vertices by assumption we have a single equation. What is even more important, knowing its solution we can reconstruct that of the original problem using the formula
 % -------------- %
\begin{eqnarray*}
\lefteqn{
\Big(\!\!
    \begin{array}{c}
    \psi(x,k)\\[-.2em] \varphi(x,k)
    \end{array}
\!\!\Big) =
\mathrm{e}^{\mp iA_j(x-j\pi)}
\bigg[
\psi_j(k)\cos k(x-j\pi)} \\ &&
+ (\psi_{j+1}( k)\mathrm{e}^{\pm iA_j\pi}-\psi_j(k)\cos k\pi)
\frac{\sin k(x-j\pi)}{\sin k\pi}
\bigg]
\,,\;\; x\in\big(j\pi, (j+1)\pi\big)\,,
\end{eqnarray*}
 % -------------- %
and in addition, the resulting function belongs to $L^p(\Gamma)$ if and only if the sequence $\{\psi_j(k)\}_{j\in\mathbb{Z}}\in \ell^p(\mathbb{Z})$,  $\:p\in\{2,\infty\}$. This relates weak solutions of the two problems but one can also make a stronger claim \cite{Pa12}:

 % -------------- %
\begin{theorem}
For any interval disjoint with the Dirichlet spectrum, \mbox{$J \subset \mathbb{R}\setminus \sigma_D$}, the restriction $(-\Delta_{\alpha,A})_J$ referring to the spectral projection $\chi_J(-\Delta_{\alpha,A})$ is unitarily equivalent to the pre-image $\eta^{(-1)}\big((L_{A})_{\eta(J)}\big)$, where $L_A$ is the operator on $\ell^2(\mathbb{Z})$ acting as $(L_A \varphi)_j = 2\cos(A_j \pi)\varphi_{j+1} + 2\cos(A_{j-1} \pi) \varphi_{j-1}$ and $\eta$ is given by~\eqref{eta}.
\end{theorem}
 % -------------- %

To proceed, let us recall another well-known mathematical problem, the \emph{almost Mathieu equation},
 % -------------- %
\begin{equation} \label{almostM}
u_{n+1} + u_{n-1} + \lambda \cos(2\pi\mu n+\theta)) u_n = \epsilon u_n\,,
\end{equation}
 % -------------- %
in particular, its critical case, $\lambda=2$, which is also referred to frequently as \emph{Harper equation}. The spectrum of the corresponding difference operator $H_{\mu,2,\theta}$ is independent of $\theta$, and plotted as a function of $\mu$ has the form of the well-known \emph{Hofstadter butterfly} \cite{Ho76} shown in Fig.~6.
 % -------------- %
\begin{figure}[h!]
\centering
    \includegraphics[scale=0.2]{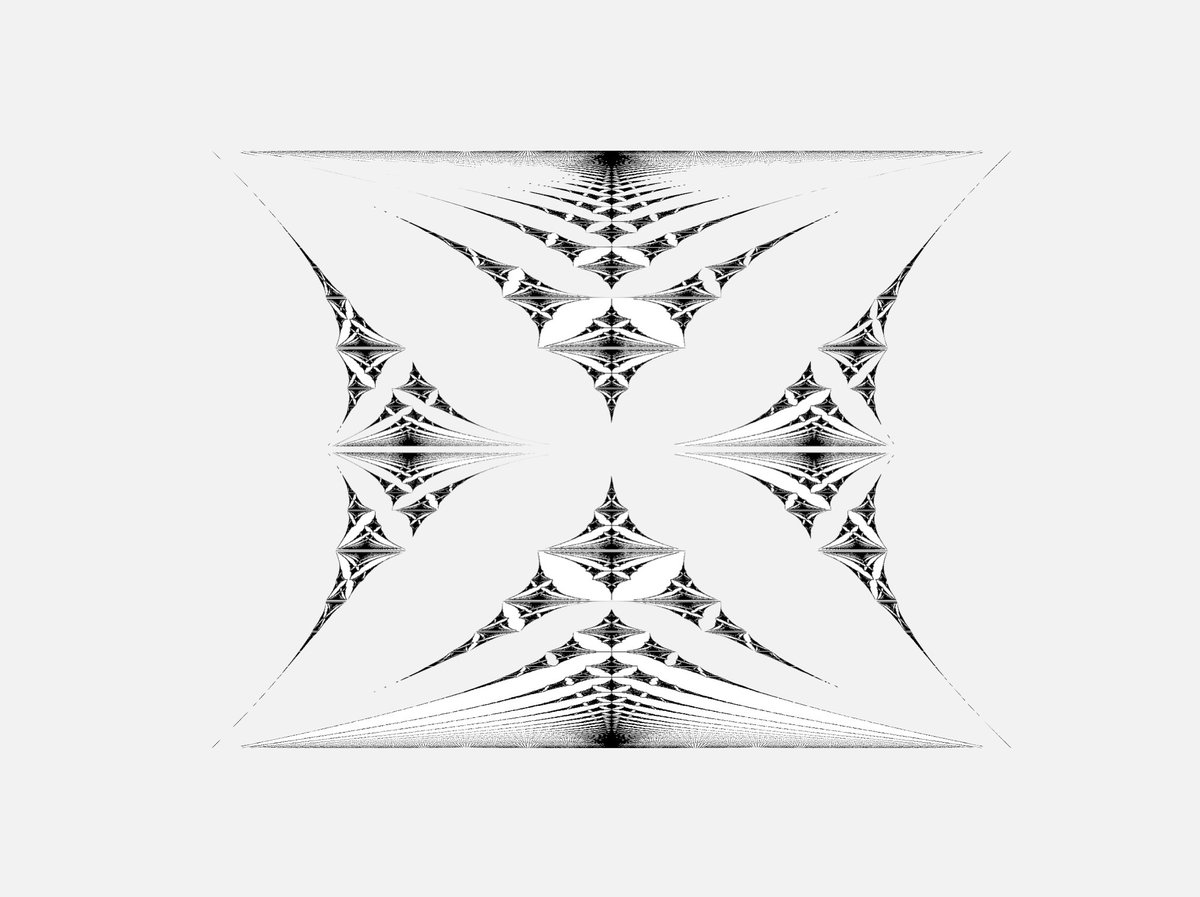}
    \caption{Hofstadter butterfly; source: Fermat's Library}
\end{figure}
 % -------------- %

The nature of the spectrum depends on the parameter involved. If $\mu\in\mathbb{Q}$, the spectrum of $H_{\mu,2,\theta}$ is easily seen to be absolutely continuous and of the band-gap type. The case of an \emph{irrational $\mu$} is much harder. The \emph{Cantor structure} was conjectured -- under the name `Ten Martini Problem' proposed by B.~Simon -- but it took two decades to achieve an affirmative answer \cite{AJ09}:

 % -------------- %
\begin{theorem}
For any $\mu\not\in\mathbb{Q}$, the spectrum of $H_{\mu,2,\theta}$ does not depend on $\theta$ and it is a Cantor
set of Lebesgue measure zero.
\end{theorem}
 % -------------- %

Let us note that the phenomenon owes its existence due to the presence of two incommensurate length scales in the problem. Such a behavior was predicted for electrons on a lattice exposed to a magnetic field by Azbel more than half a century ago \cite{Az64}; recently it was observed experimentally by several groups on graphene lattices \cite{De+13, Po+13}.

The question is how is this related to our problem. To see that we employ a trick originally proposed by M.~Shubin \cite{Sh94}: we consider a rotation algebra $A_\mu$ generated by elements $u,v$ such that $uv= \mathrm{e}^{2\pi i\mu}vu$. Such an algebra is simple for $\mu\not\in\mathbb{Q}$, thus having faithful representations. We construct two representations of $A_\mu$ which map a single element of the algebra, namely $u+v+u^{-1}+v^{-1}\in A_\mu$, to operators $L_A$ and $H_{\mu,2,\theta}$, respectively, and use it to conclude that their spectra coincide, $\sigma(L_A)=\sigma(H_{\mu,2,\theta})$.
This allows us using the above described duality and the fact that the function \eqref{eta} is locally analytic to obtain, in a cheap way indeed, the following nontrivial result \cite{EV17}:

 % -------------- %
\begin{theorem}
If $\alpha \ne 0$ and $\mu \notin \mathbb{Q}$, then $\sigma(-\Delta_{\alpha,A})$ does not depend on $\theta$ and it is a disjoint union of the isolated-point family $\{n^2 |\, n\in\mathbb{N}\}$ and \emph{Cantor sets}, one inside each interval $(-\infty, 1)$ and $\big(n^2, (n+1)^2\big)$, $n \in \mathbb{N}$. Moreover, the overall Lebesgue measure of $\sigma(-\Delta_{\alpha,A})$ is zero.
\end{theorem}
 % -------------- %

We can also translate other results concerning the critical almost Mathieu operator to the problem of a chain graph with the linear magnetic field. For instance, we have \cite{LS16}:

  % -------------- %
\begin{theorem}
Let $A_j = \mu j + \theta$ hold for some $\mu, \theta \in \mathbb{R}$ and every $j\in\mathbb{Z}$. There exists a dense $G_\delta$ set of the slopes $\mu$ for which, and for all $\theta$, the Hausdorff dimension $\dim_H \sigma(-\Delta_{\alpha,A}) = 0$.
\end{theorem}
 % -------------- %

On the other hand, using the result \cite{HLQZ19} we arrive at the following claim:

 % -------------- %
\begin{theorem}
There is \emph{another} dense set of the slopes $\mu$, with positive Hausdorff measure, for which, on the contrary,  $\dim_H \sigma(-\Delta_{\alpha,A}) > 0$.
\end{theorem}
 % -------------- %
Taken together, the last two results show in a convincing way how subtle the spectral problem in this system is.

%%%%%%%%%%%%%%%%%%%%%%%%%%%%%%%%%%%%%%%%%%%%%%
\section{Bethe-Sommerfeld graphs} \label{s:BS}
\setcounter{equation}{0}

Let us now change the topic. Figuratively speaking, the graphs in the previous example had `many' gaps indeed, now we will deal in contrast with periodic graphs having `just a few' gaps. We mentioned in the introduction the \emph{Bethe--Sommerfeld conjecture}, made in early days of quantum mechanics and  stating that systems periodic in more than one direction have at most a finite number of open spectral gaps. The reasoning behind the conjecture is simple: in one dimension the dispersion curves of the free system do not have many intersections, and therefore it is not difficult for a perturbation to open gaps even at high energies, in sheer contrast to higher dimensional systems. However, to establish such a claim mathematically proved to be a hard problem and the first rigorous results appeared only in the 1980s. Since then a substantial progress was made and nowadays the validity of the conjecture is established for wide classes of operators \cite{Pa08}.

It is natural to ask whether the same could be true for quantum graphs which one can regard as systems mixing different dimensionalities. The common understanding is that, while the above described heuristic argument applies here too, the finiteness of the gap number `is not a strict law' \cite{BK13}. The reason comes again from the topology: one is able to create an infinite number of gaps in the spectrum of a periodic graph by \emph{decorating} its vertices by copies of a fixed compact graph. This fact was observed first in the combinatorial graph context \cite{SA00} and the argument extends easily to metric graphs we consider here \cite{Ku04} in the way illustrated in Fig.~7.
 % -------------- %
\begin{figure}[h!]
\centering
    \includegraphics[scale=0.45,angle=0]{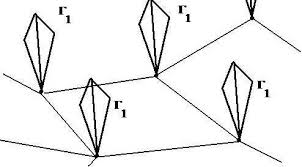}
    \caption{Decoration of a graph; courtesy: Peter Kuchment}
\end{figure}
 % -------------- %

Thus the question rather could be whether it is a `law' at all: do infinite periodic graphs having a finite and nonzero number of open gaps exist? From obvious reasons we would call them \emph{Bethe-Sommerfeld graphs}.

The answer depends on the vertex coupling. Recall that the general self-adjoint matching conditions are described by relation \eqref{bc}. One can decompose them into the \emph{Dirichlet}, \emph{Neumann}, and \emph{Robin} parts \cite[Thm.~1.4.4]{BK13} corresponding to the eigenspaces of $U$ with eigenvalues $-1,\,1$, and the rest, respectively. If the latter is absent we call such a coupling \emph{scale-invariant}; as an example, one can mention the Kirchhoff coupling. Graphs with such vertices do not have the property we seek \cite{ETu18}:
 % -------------- %
\begin{theorem}
An infinite periodic quantum graph does \emph{not}  belong to the Bethe- Sommerfeld class if the couplings at its vertices are scale-invariant.
\end{theorem}
 % -------------- %
Furthermore, recall the result of \cite{BB13} which shows in a particularly illuminating way that the probability of finding an energy value in a band is universal within wide classes of periodic graphs. If this probability is less than one, such a graph clearly has an infinite number of gaps, in the opposite case one might be tempted to believe that there are no gaps.

Nevertheless, the answer to our question is affirmative \cite{ETu18}:
 % -------------- %
\begin{theorem} \label{thm:BSex}
Bethe--Sommerfeld graphs exist.
\end{theorem}
 % -------------- %
Since it is an existence claim it is sufficient, of course, to demonstrate an example. With this aim we revisit the model of a \emph{rectangular lattice graph} with a $\delta$ coupling in the vertices, sketched in Fig.~8, studied in \cite{Ex96b, EG96}.
 % -------------- %
\begin{figure}[h!]
\centering
    \includegraphics[trim=4cm 23.5cm 14cm 2.5cm, width=3cm]{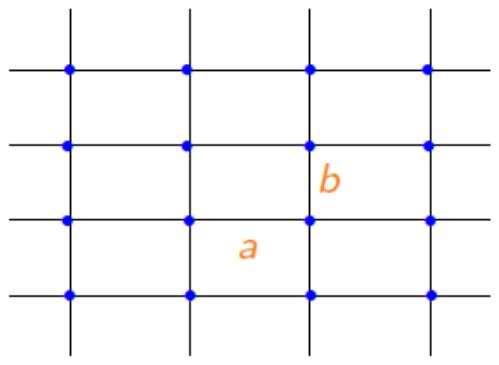}
    \caption{Rectangular lattice graph}
\end{figure}
 % -------------- %
A number $k^2>0$ belongs to a gap if and only if $k>0$ satisfies the \emph{gap condition} which is easily derived; it reads
 % -------------- %
\begin{equation} \label{gap+}
2k\bigg[\tan\left(\frac{ka}{2}-\frac{\pi}{2}\left\lfloor\frac{ka}{\pi}\right\rfloor\right)
+\tan\left(\frac{kb}{2}-\frac{\pi}{2}\left\lfloor\frac{kb}{\pi}\right\rfloor\right) \bigg]<\alpha \;\quad \text{for}\;\; \alpha>0
\end{equation}
 % -------------- %
and
 % -------------- %
\begin{equation} \label{gap-}
2k\bigg[\cot\left(\frac{ka}{2}-\frac{\pi}{2}\left\lfloor\frac{ka}{\pi}\right\rfloor\right)
+\cot\left(\frac{kb}{2}-\frac{\pi}{2}\left\lfloor\frac{kb}{\pi}\right\rfloor\right) \bigg]<|\alpha| \;\quad \text{for}\;\; \alpha<0\,;
\end{equation}
 % -------------- %
we neglect the Kirchhoff case, $\alpha=0$, where $\sigma(H)=[0,\infty)$. We note that for $\alpha<0$ the spectrum extends to the negative part of the real axis and may have a gap there, which is not important here because there is not more than a single negative gap, and this gap always extends to positive values.

Let us recall briefly what is known about the spectrum of the corresponding Hamiltonian $H\,$ \cite{Ex96b, EG96}: it depends on the ratio $\theta=\frac{a}{b}$. To begin with, $\sigma(H)$ has clearly infinitely many gaps if $\theta$ is rational, except for the Kirchhoff case, $\alpha=0$, when $\sigma(H)=[0,\infty)$. The same is true if $\theta$ is is \emph{an irrational well approximable by rationals}, which means that in the continued fraction representation, $\theta=[a_0;a_1,a_2,\dots]$, the coefficient sequence $\{a_j\}$ is unbounded.

On the other hand, $\theta\in\mathbb{R}$ is said to be \emph{badly approximable} if the sequence is bounded; according to a well know number theory result this fact is equivalent to the existence of a $c>0$ such that
 % -------------- %
$$
\Big|\theta-\frac{p}{q}\Big|>\frac{c}{q^2}
$$
 % -------------- %
for all $p,q\in\mathbb{Z}$ with $q\neq0$. For such numbers we define \emph{Markov constant} by
 % -------------- %
\begin{equation} \label{Markov}
\mu(\theta):=
\inf\bigg\{c>0\;\big|\;\left(\exists_\infty(p,q)\in\mathbb{N}^2\right) \Big(\Big|\theta-\frac{p}{q}\Big|<\frac{c}{q^2}\Big)\bigg\}\,,
\end{equation}
 % -------------- %
with $\exists_\infty$ meaning `there exist infinitely many', noting that $\mu(\theta)=\mu(\theta^{-1})$, and its `one-sided analogues' referring to similar approximations with the positive and negative values of $\theta-\frac{p}{q}$, respectively.

As an example, consider the \emph{golden mean}, $\theta=\frac{\sqrt{5}+1}{2} = [1;1,1,\dots]$, which can be regarded as the `worst' irrational. A quick look reveals interesting properties \cite{EG96} which follow from the behavior of the left-hand side of \eqref{gap+}. It is clear from the formula that, as a function of $k$, it has a saw-tooth shape with a number of local minima. In view of \eqref{gap+} it is the latter which determine the existence of gaps. We plot them in Fig.~9; since the occurrence of small values of those minimum points becomes rare as the energy grows we use the logarithmic scale of $k$ in the graph.
 % -------------- %
\begin{figure}[h!]
\centering
    \includegraphics[scale=0.65,angle=0]{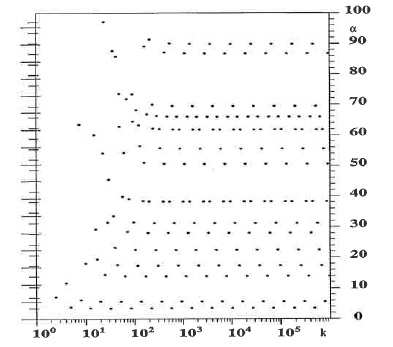}
    \caption{Solution of the condition \eqref{gap+}: spectral gaps are situated around the values of $k$ for which the plotted minima are smaller than $\alpha$.}
\end{figure}
 % -------------- %

For comparison it is useful to note that for a rational $\theta$ the left-hand side of \eqref{gap+} would have a periodic sequence of zeros giving rise to gaps for any $\alpha>0$, and furthermore, for an irrational $\theta$ that is well approximated by rationals there will be no zeros but one could find a sequence of the local minima tending to zero, and consequently, one would have again an infinite series of gaps for any $\alpha>0$.

In contrast, in our case the set of the minima has positive lower bound, which means that there are no gaps in the spectrum for $\alpha$ small enough. In addition, the picture shows that making the value of the coupling constant larger, one can open other infinite gap series: they appear at the values $\alpha=\frac{\pi^2}{\sqrt{5ab}} \theta^{\pm1/2} |n^2-m^2-nm|,\: n,m\in\mathbb{N}\:$. However, the lowest series of the function minimum values approaches the limit as $k\to\infty$ \emph{from above} so no situation with a finite number of gaps arises.

But a closer look shows a more complex picture \cite{ETu18}:
 % -------------- %
\begin{theorem}
Let $\frac{a}{b}=\theta=\frac{\sqrt{5}+1}{2}$, then the following claims are valid:
 % -------------- %
\begin{itemize}
\item [(i)] If $\alpha>\frac{\pi^2}{\sqrt{5}a}$ or $\alpha\leq-\frac{\pi^2}{\sqrt{5}a}$, there are infinitely many spectral gaps.
\item [(ii)] If
 % -------------- %
$$
-\frac{2\pi}{a}\tan\Big(\frac{3-\sqrt{5}}{4}\pi\Big)\leq\alpha\leq\frac{\pi^2}{\sqrt{5}a}\,,
$$
 % -------------- %
there are no gaps in the positive spectrum.
\item[(iii)] If
 % -------------- %
$$
-\frac{\pi^2}{\sqrt{5}a}<\alpha<-\frac{2\pi}{a}\tan\Big(\frac{3-\sqrt{5}}{4}\pi\Big),
$$
 % -------------- %
there is \emph{a nonzero and finite number of gaps} in the positive spectrum.
\end{itemize}
 % -------------- %
\end{theorem}
 % -------------- %

As a corollary, this result establishes the validity of Theorem~\ref{thm:BSex}. The window in which the golden-mean lattice has the Bethe-Sommerfeld property is narrow, corresponding roughly to the values $4.298 \lesssim -\alpha a \lesssim 4.414$. Within these limits, we can do even better and  control the number of gaps in the Bethe-Sommerfeld regime \cite{ETu18}:
 % -------------- %
\begin{theorem}
For a given $N\in\mathbb{N}$, there are \emph{exactly $N$ gaps} in the positive spectrum if and only if $\alpha$ satisfies $-A_{N+1}\le \alpha <A_N$, where
 % -------------- %
$$
A_{j}:=\frac{2\pi\left(\theta^{2j}-\theta^{-2j}\right)}{\sqrt{5}}\tan\left(\frac{\pi}{2}\theta^{-2j}\right)\,, \quad j\in\mathbb{N}\,.
$$
 % -------------- %
\end{theorem}
 % -------------- %
Note that the numbers $A_j$ form an increasing sequence the first element of which is $A_1=2\pi\tan\left(\frac{3-\sqrt{5}}{4}\pi\right)$ and $A_j<\frac{\pi^2}{\sqrt{5}}$ holds for all $j\in\mathbb{N}$.

Proofs of the above results are based on properties of Diophantine approximations. In a similar way one can get a more general result \cite{ETu18}:

 % -------------- %
\begin{theorem}
Let $\theta=\frac{a}{b}$ and define

\vspace{-1em}

{\scriptsize

 % -------------- %
$$
\gamma_+:=\min\left\{\inf_{m\in\mathbb{N}}\left\{\frac{2m\pi}{a}\tan\left(\frac{\pi}{2}(m\theta^{-1}-\lfloor m\theta^{-1}\rfloor)\right)\right\},\inf_{m\in\mathbb{N}}\left\{\frac{2m\pi}{b}\tan\left(\frac{\pi}{2}(m\theta-\lfloor m\theta\rfloor)\right)\right\}\right\}
$$
 % -------------- %

}

and $\gamma_-$ similarly with $\lfloor\cdot\rfloor$ replaced by $\lceil\cdot\rceil$ and $\tan$ by $-\tan$. If $\pm\alpha>0$ and the coupling constant satisfies the inequality
 % -------------- %
$$
\gamma_\pm<\pm\alpha<\frac{\pi^2}{\max\{a,b\}}\mu(\theta)\,,
$$
 % -------------- %
then there is \emph{a nonzero and finite number of gaps} in the positive spectrum.
\end{theorem}
 % -------------- %
This allows us to construct further examples, in particular, to show that also lattices with \emph{repulsive $\delta$ coupling}, $\alpha>0$, may exhibit the Bethe-Sommerfeld property, cf.~\cite{ETu18} for more details.

%%%%%%%%%%%%%%%%%%%%%%%%%%%%%%%%%%%%%%%%%%%%%%%%%%%%%%%%%%%%%%%%%%%%%%%%%%%%%%%%%%
\section{An interlude about the meaning of the vertex coupling} \label{s:coupling}
\setcounter{equation}{0}

Our next aim is to show one more example where a topological characteristics, in this case the vertex degree, or rather its \emph{parity}, has a substantial influence on the spectrum. Before coming to it we have to look more closely to the meaning of the vertex coupling \eqref{bc} which contain a number of parameters. One way to determine which are the `right' ones was suggested longtime ago, following the original proposal of Linus Pauling to use metric graphs as an instrument to analyze molecules of aromatic hydrocarbons \cite{Pa36}. Since such graphs are supposed to model thin networks, it seems natural to inspect how the Laplacian, or more generally a Schr\"odinger operator, supported by such a network behaves when the width of the network threads tends to zero. For the Neumann Laplacian a formal use of the Green's formula suggests that such a squeezing limit would lead to the \emph{Kirhhoff} coupling \cite{RS53}. This guess proved to be correct but it took time before the claim was established rigorously \cite{EP05, FW93, KZ01, RS01, Sa00}.

For other couplings such a straightforward procedure based on the Laplacian is not sufficient. A modification leading to the \emph{$\delta$ coupling} is easy: one has to employ instead the Schr\"odinger operator with a family of \emph{scaled potentials}, the integral of which yields the coupling constant \cite{EP09}. The general case is considerably more complicated: one has not only to add more potentials as well as \emph{magnetic fields}, but also to modify locally the \emph{network topology}. One employs a family of magnetic Schr\"odinger operator supported by a \emph{Neumann-type} network sketched in Fig.~10 with \emph{Neumann-type boundary}. This includes both a network-type subset of $\mathbb{R}^3$ supporting a Neumann Laplacian as well as the situation when the manifold in question is the surface of such a region having thus no boundary.
 % -------------- %
\begin{figure}[h!]
\centering
    \includegraphics[trim=4cm 21.5cm 5cm 2.5cm, width=12cm]{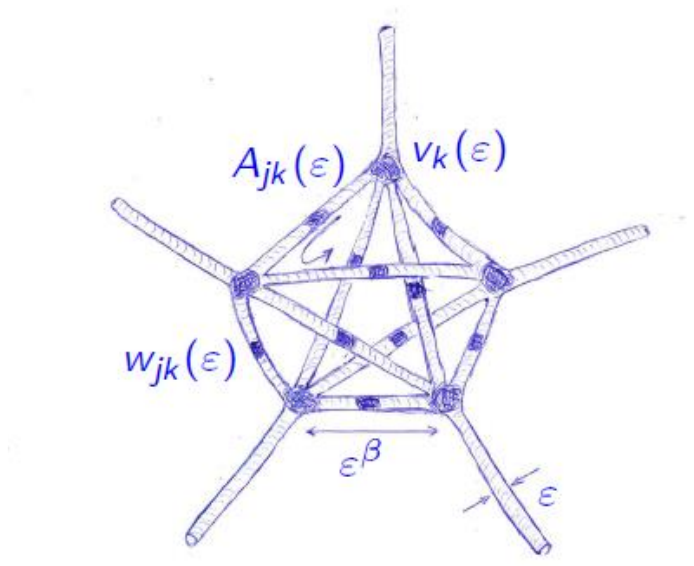}
    \caption{Squeezed network approximation: the manifold supports scalar potentials $v_k$ in a neighborhood of the endpoints of the connecting tubes and $w_{jk}$ in the vicinity of their midpoints, together with the vector potentials $A_{jk}$, all dependent on the tube radius $\varepsilon$; the length of the connecting tubes shrinks as $\varepsilon^\beta$ when $\varepsilon\to 0$.}
\end{figure}
 % -------------- %
Choosing properly the potentials, both the scalar supported by the shaded regions and the vector ones, as functions of $\varepsilon$, singular as $\varepsilon\to 0$, and assuming that the network width shrinks much faster than the length of the added edges, $\beta<\frac{1}{13}$, one can approximate \emph{any} vertex coupling in the \emph{norm-resolvent sense} \cite{EP13}.

Let us note in passing that the analogous squeezing problem for a network-shaped region supporting the Laplacian with other boundary conditions, in particular, the Dirichlet one, is much harder. The ways to get a nontrivial limit choosing a proper energy renormalization are known \cite{CE07, Gr08, MV07}, however, we are far from a complete understanding similar to the Neumann case. We are not going to discuss this problem here.

One has to acknowledge that the above described construction is rather a baroque one. Alternatively, one can take an pragmatic approach and look whether a given vertex coupling could be useful to model a particular system. To motivate the problem discussed in the next section let us mention one situation where a quantum graph model was used recently \cite{SK15}. It concerns the well-known \emph{Hall effect} sketched in Fig.~11.
 % -------------- %
\begin{figure}[h!]
\centering
    \includegraphics[width=0.45\textwidth]{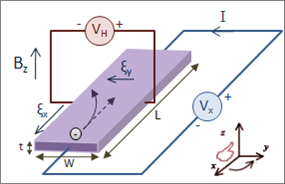}
    \caption{The Hall effect; source: Wikipedia}
\end{figure}
 % -------------- %
It was discovered 140 years ago when Edwin Hall noticed that magnetic field induces a \emph{voltage perpendicular} to the current. However, the true excitement came in the 1980s when Klaus von Klitzing and others discovered that in the quantum regime the Ohm's law ceased to hold and the corresponding  conductivity was \emph{quantized} with a great precision -- this fact lead already to two Nobel Prizes. What is even more intriguing, however, is that in ferromagnetic materials one can also observe a similar behavior in the \emph{absence of external magnetic field} -- this is referred to as the \emph{anomalous Hall effect}.

The mechanism of the `usual' quantum Hall effect is well understood -- this was the second one of the mentioned Nobel prizes -- which cannot be said about its anomalous counterpart of which it is conjectured that it comes from internal magnetization in combination with the spin-orbit interaction.
Recently a model was proposed in which the material structure of the sample is described by lattice of $\delta$-coupled rings, sketched in Fig.~12,
 % -------------- %
\begin{figure}[h!]
\centering
    \includegraphics[clip, trim=6.5cm 3cm 7cm 3cm,angle=270,width=1\textwidth]{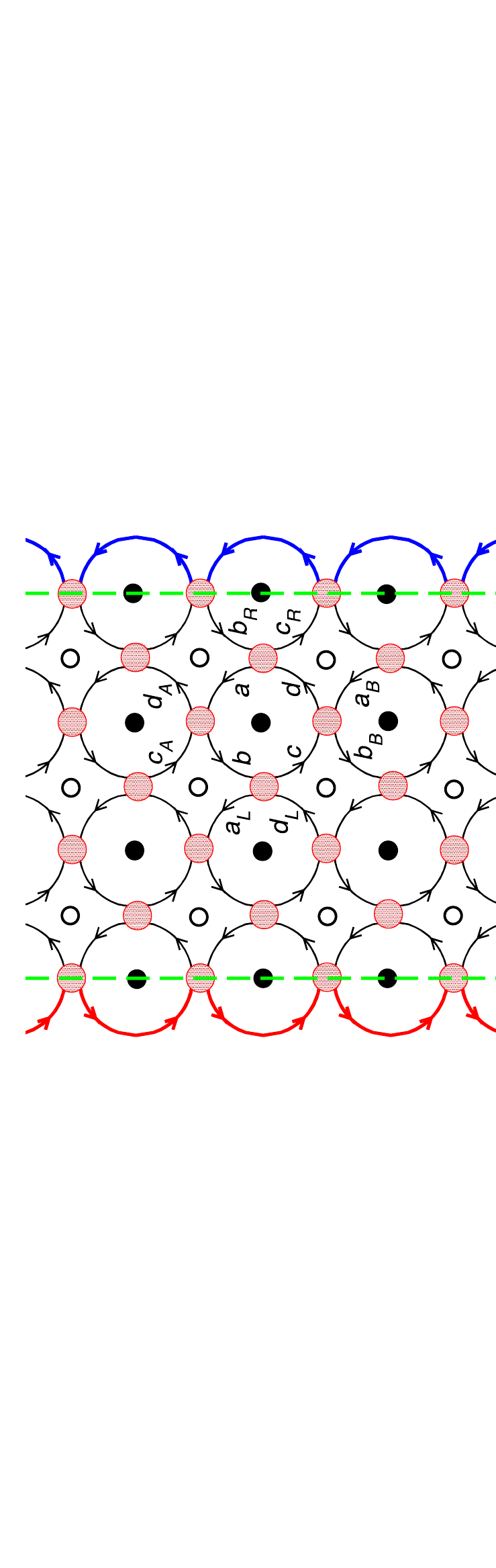}
    \caption{A model of anomalous Hall effect; source \cite{SK15}}
\end{figure}
 % -------------- %
topologically equivalent to a square lattice. Looking at the picture we recognize a problem: to mimick the rotational motion of atomic orbitals responsible for the magnetization, the authors had to impose `by hand' the requirement that the electrons move only one way, anticlockwise, on the loops of the lattice. Such an assumption, however, cannot be justified from the first principles. On the other hand, it \emph{is} possible to break the time-reversal invariance, not at graph edges but at its \emph{vertices}; the next section is devoted to the discussion of a model with this property.

%%%%%%%%%%%%%%%%%%%%%%%%%%%%%%%%%%%%%%%%%%%%%%%%%%%%%%%%%%%%%%%%%%%%%%%
\section{Violating the time-reversal invariance} \label{s:timereversal}
\setcounter{equation}{0}

To construct an example of vertex coupling with the indicated property, recall that on a star graph consisting on $N$ halfline meeting at a single vertex with the coupling \eqref{bc} the \emph{on-shell S-matrix}, that is, the component of the direct-integral decomposition of the scattering operator referring to the momentum $k\:$ \cite[Sec.~2.1]{BK13} equals
 % ------------- %
\begin{equation} \label{onshell}
S(k) = \frac{k-1 +(k+1)U}{k+1 +(k-1)U}\,,
\end{equation}
 % ------------- %
in particular, we have $U=S(1)$. The maximum preferential orientation of such a vertex at $k=1$ would be thus achieved if we choose
{\small
 % ------------- %
\begin{equation} \label{rotU}
U = \left( \begin{array}{ccccccc}
0 & 1 & 0 & 0 & \cdots & 0 & 0 \\ 0 & 0 & 1 & 0 & \cdots & 0 & 0 \\ 0 & 0 & 0 & 1 & \cdots & 0 & 0 \\ \cdots & \cdots & \cdots & \cdots & \cdots & \cdots & \cdots \\ 0 & 0 & 0 & 0 & \cdots & 0 & 1 \\ 1 & 0 & 0 & 0 & \cdots & 0 & 0
\end{array} \right).
\end{equation}
 % ------------- %
}
To find properties of the vertex coupling \eqref{bc} with the matrix \eqref{rotU} we consider again the star graph. Writing the coupling conditions componentwise, we have
 % ------------- %
\begin{equation} \label{rotat}
(\psi_{j+1}-\psi_j) + i(\psi'_{j+1}+\psi'_j) = 0\,, \quad j\in\mathbb{Z}\; (\mathrm{mod}\,N)\,,
\end{equation}
 % ------------- %
which is non-trivial for $N\ge 3$ and obviously non-invariant with respect to the reverse in the edge numbering order, or equivalently, with respect to the complex conjugation representing the time reversal. For such a star-graph Hamiltonian we obviously have $\sigma_\mathrm{ess}(H) = \mathbb{R}_+$. It is also easy to check that $H$ has eigenvalues $-\kappa^2$, where
 % ------------- %
\begin{equation} \label{star-ev}
\kappa= \tan \frac{\pi m}{N}
\end{equation}
 % ------------- %
with $m$ running through $1,\dots,[\frac{N}{2}]$ for $N$ odd and $1,\dots,[\frac{N-1}{2}]$ for $N$ even. Thus $\sigma_\mathrm{disc}(H)$ is \emph{always nonempty}, in particular, $H$ has a single negative eigenvalue for $N=3,4$ which is equal to $-3$ and $-1$, respectively.

Let us turn to the scattering. Looking at the relation \eqref{onshell} it might seem that transport becomes trivial at small and high energies, $\lim_{k\to 0} S(k)=-I$ and $\lim_{k\to\infty} S(k)=I$. However, caution is needed; the formal limits lead to a false result if $+1$ or $-1$ are eigenvalues of $U$. A simple counterexample is the (scale invariant) Kirchhoff coupling for which $U$ has only $\pm 1$ as its eigenvalues; the on-shell S-matrix is then independent of~$k$ and it is \emph{not} a multiple of the identity. A straightforward computation yields the explicit form of $S(k)$: denoting for simplicity
 % ------------- %
$
\eta := \frac{1-k}{1+k}
$
 % ------------- %
we have \cite{ETa18}
 % ------------- %
$$
S_{ij}(k) = \frac{1-\eta^2}{1-\eta^N} \left\{ -\eta\, \frac{1-\eta^{N-2}}{1-\eta^2}\,\delta_{ij} +
(1-\delta_{ij})\, \eta^{(j-i-1)(\mathrm{mod}\,N)} \right\}.
$$
 % ------------- %
This formula suggests, in particular, that the S-matrix high-energy behavior, $\eta\to -1$, could be determined by the \emph{parity} of the vertex degree $N$. In the cases with the lowest $N$ we get
{\small
 % ------------- %
$$
S(k)= \frac{1+\eta}{1+\eta+\eta^2} \left( \begin{array}{ccc}
-\frac{\eta}{1+\eta} & 1 & \eta \\ \eta & -\frac{\eta}{1+\eta} & 1 \\ 1 & \eta & -\frac{\eta}{1+\eta}
\end{array} \right)
$$
 % ------------- %
}
and
{\small
 % ------------- %
$$
S(k)= \frac{1}{1+\eta^2} \left( \begin{array}{cccc}
-\eta & 1 & \eta & \eta^2 \\ \eta^2 & -\eta & 1 & \eta \\ \eta & \eta^2 & -\eta & 1 \\ 1 & \eta & \eta^2 & -\eta
\end{array} \right)
$$
 % ------------- %
}
for $N=3,4$, respectively. We see that $\lim_{k\to\infty} S(k)=I$ holds for $N=3$ and more generally \emph{for all odd $N$}, while for the \emph{even ones} the limit is \emph{not a multiple of identity}; this comes from the fact that in the latter case $U$ has both $\pm 1$ as its eigenvalues, while for $N$ odd $-1$ is missing.

Let us look how this property influences spectra of periodic quantum graphs. To this aim we compare two lattices shown in Fig.~13, the square one and the honeycomb -- or hexagonal -- one, with their period cells.
 % ------------- %
\begin{center}
\begin{figure}
\hspace{-4em}\includegraphics[clip, trim=3.5cm 11.5cm 7cm 11.5cm,width=.6\textwidth]{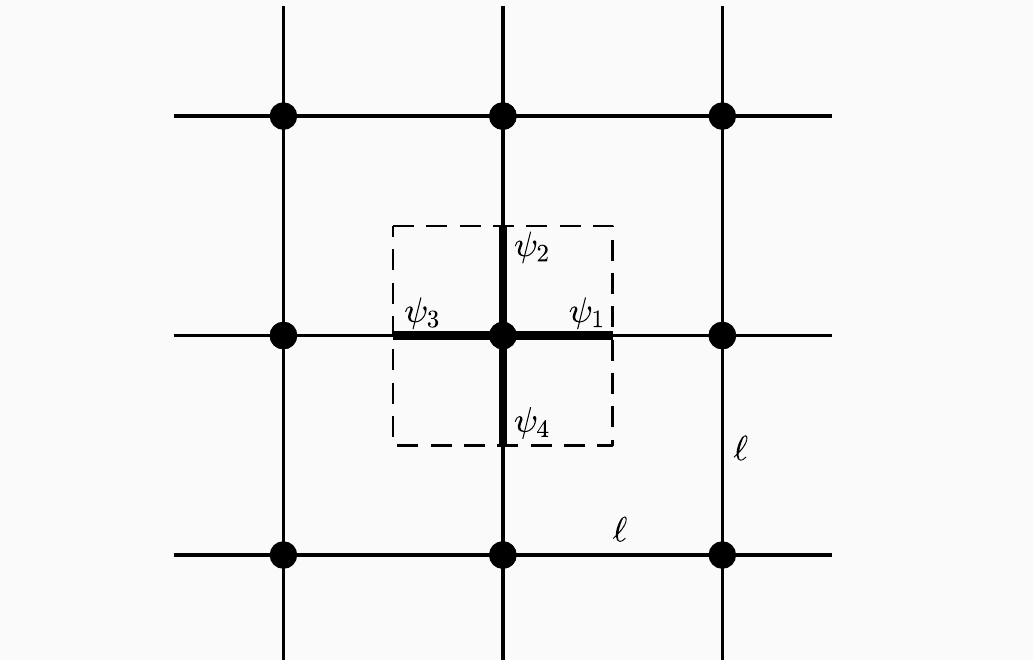}\hspace{-2em}\includegraphics[clip, trim=1.5cm 15.9cm 10.6cm 6.5cm, width=0.45\textwidth]{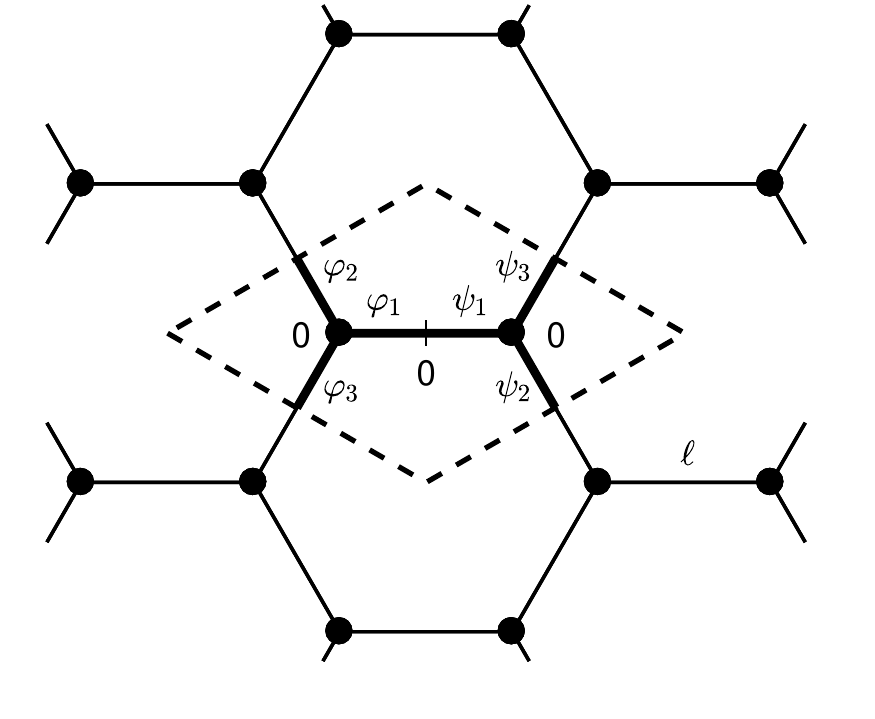}
\caption{The square and honeycomb lattices}
\end{figure}
\end{center}
 % ------------- %

\vspace{-2.5em}

\noindent The secular equations determining the spectrum for the two cases are easy to derive,
 % ------------- %
\begin{equation} \label{sec-square}
16i\,\mathrm{e}^{i(\theta_1+\theta_2)}\, k\, \sin k\ell \big[ (k^2-1) (\cos\theta_1 + \cos\theta_2) + 2(k^2+1) \cos k\ell \big]=0
\end{equation}
 % ------------- %
and respectively
{\small
 % ------------- %
\begin{equation} \label{sec-honey}
16i\,\mathrm{e}^{-i(\theta_1+\theta_2}\,k^2\sin k\ell\, \Big( 3 + 6k^2 - k^4 +4d_\theta(k^2-1) + (k^2+3)^2 \cos 2k\ell \Big) = 0\,,
\end{equation}
 % ------------- %
}
where
 % ------------- %
$
d_\theta := \cos\theta_1 + \cos(\theta_1-\theta_2) + \cos\theta_2
$
 % ------------- %
and $\frac{1}{\ell}(\theta_1,\theta_2) \in [-\frac{\pi}{\ell},\frac{\pi}{\ell}]^2$ is the quasimomentum. They are tedious to solve except the \emph{flat band cases}, $\sin k\ell=0$, however, we can present the solutions to \eqref{sec-square} and \eqref{sec-honey} in a \emph{graphical form} as shown in Figures 14 and 15 for two values of the lattice spacing, namely $\ell=\frac32$ (solid curve) and $\ell=\frac14$ (dashed one).
 % -------------- %
\begin{figure}[htbp]
     \centering
     \includegraphics[clip, trim=1.5cm 6.8cm 4cm 15cm,width=0.75\textwidth]{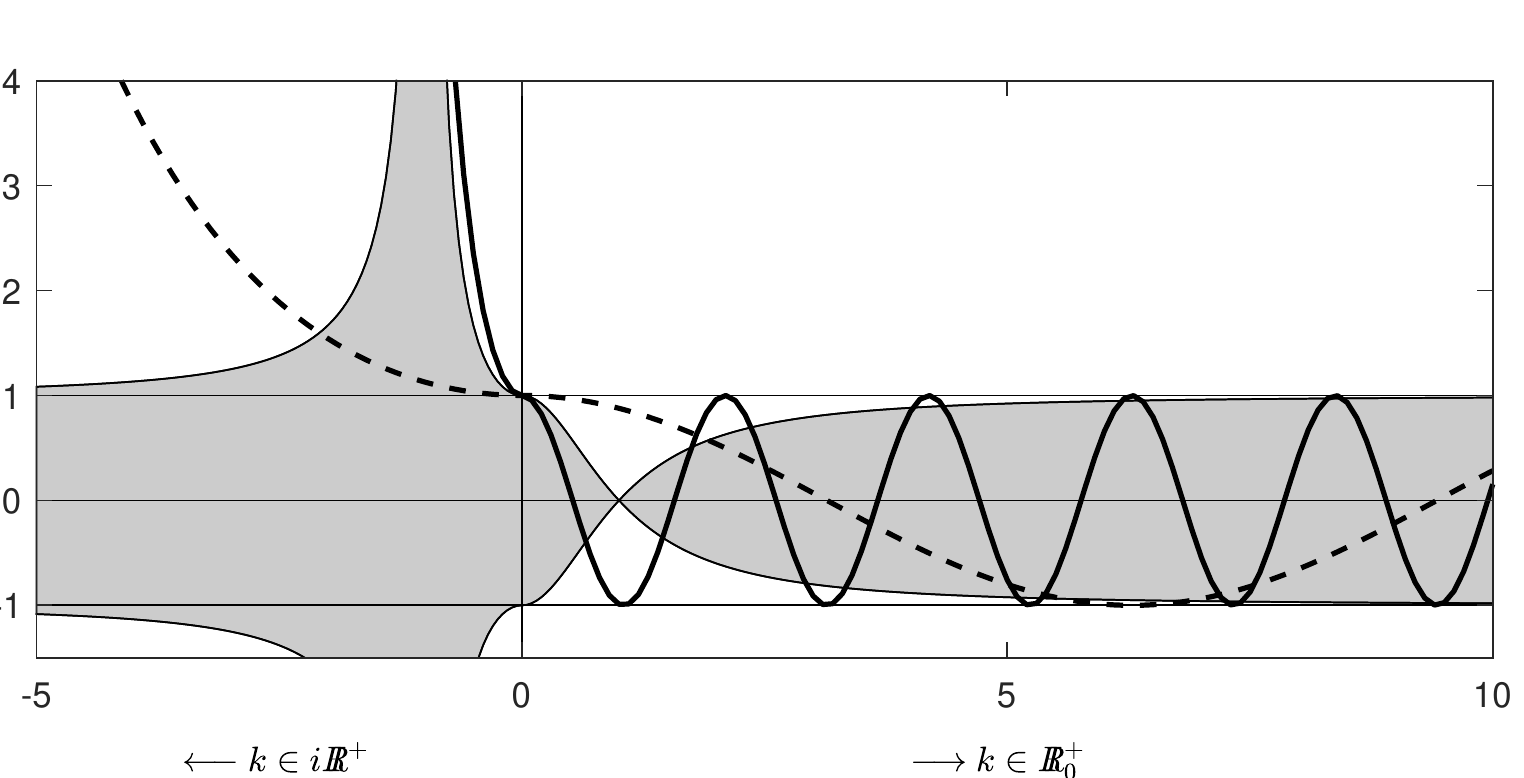}
     \caption{Graphical solution of the square lattice spectrum}
\end{figure}
 % -------------- %
To be specific, to make the last factor on the left-hand sides of the conditions \eqref{sec-square} and \eqref{sec-honey} vanish at the energy $k^2$ requires $\cos jk\ell = f_j(\theta_1,\theta_2;k),\: j=1,2$, for appropriately chosen functions $f_j$; in the negative part of the spectrum we have to replace trigonometric functions by hyperbolic ones. Plotting the ranges of the functions $f_1,\,f_2$ as shaded regions in Figures 14 and 15, respectively, we find the solutions to the secular equations as the points on the $x$-axis for which cosine/hyperbolic cosine value lies in the shaded region.
 % -------------- %
\begin{figure}[htbp]
     \centering
     \includegraphics[clip, trim=1.5cm 6.8cm 4cm 15cm,width=.75\textwidth]{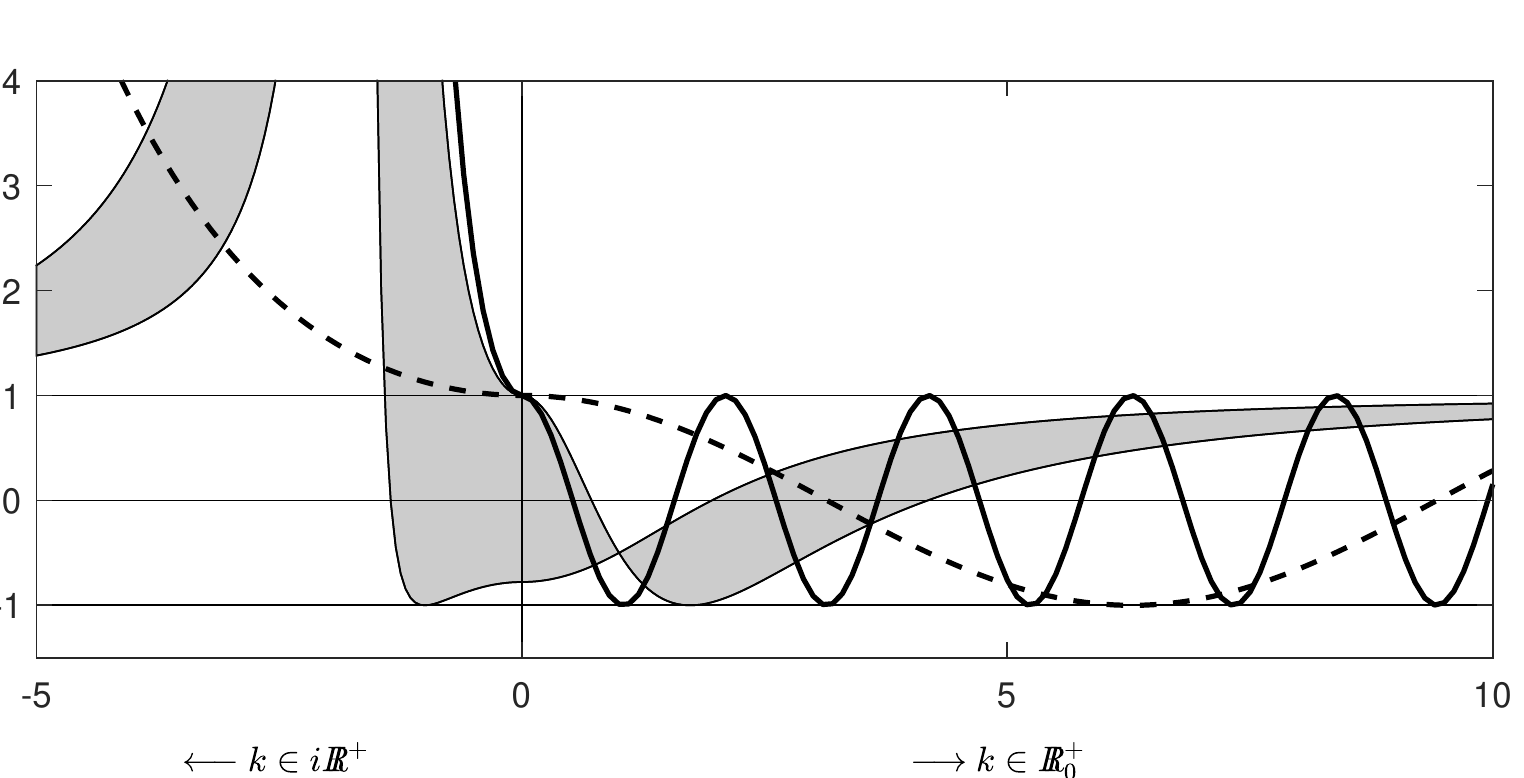}
     \caption{Graphical solution of the hexagonal lattice spectrum}
\end{figure}
 % -------------- %

Comparing the two systems we see that some features are common:

\vspace{-.5em}

\begin{itemize}
\setlength{\itemsep}{0pt}

\item the number of open gaps is \emph{always infinite},

\item the gaps occur in pairs with a flat band in the middle,

\item for some values of $\ell$ a band may \emph{degenerate},

\item the negative spectrum is \emph{always nonempty} (for large $\ell$ it is not seen from the pictures but reader would have no difficulty understanding that the graph of the hyperbolic cosine crosses the shaded area. Looking into the shape of the latter around the peak, it is straightforward to check that the gaps become \emph{exponentially narrow} around star graph eigenvalues given by \eqref{star-ev} as $\ell\to\infty$.

\end{itemize}

\vspace{-.5em}

However, the \emph{high energy behavior} of these lattices is substantially different:

\vspace{-.5em}

\begin{itemize}
\setlength{\itemsep}{0pt}

\item the spectrum is dominated by \emph{bands} for square lattices, with the gap widths (in the energy variable) tending to a nonzero constant as the gap index increases,

\item in contrast, for hexagonal lattices the spectrum is dominated by \emph{gaps} for hexagonal lattices, the bands being now asymptotically of a constant width.

\end{itemize}

\vspace{-.5em}

In the light of what we said about the properties of the on-shell S-matrix the difference is easy to be understood. Vertices of the square lattice are degree four, hence a particle of a high energy approaching a vertex has roughly the same probability to leave it in any direction. On the other hand, at a vertex of degree three as for the honeycomb lattice, the reflection is dominating so -- except for the narrow bands -- the particle remains trapped at a single edge. As for the asymptotic band behavior in the negative part of the spectrum, we note that the edges represent now, in the physicist's language, classically forbidden zones and the transport is possible only due to tunneling between the vertices which becomes more difficult as $\ell$ increases.

One can also discuss more general vertex coupling violating the time-rever\-sal invariance, for instance, an \emph{interpolation} between the $\delta$-coupling and the one considered here, but we are not going do it here referring to \cite{ETT18} for the corresponding analysis.

%%%%%%%%%%%%%%%%%%%%%%%%%%%%%%%%%%%%%%%%%%%%%%%%
\section{Band edge positions} \label{s:bandedge}
\setcounter{equation}{0}

The vertex coupling characterized by the matrix \eqref{rotU} can serve also another purpose. To introduce the topic, let us first recall that looking for edges of spectral bands of a one-dimensional periodic Schr\"odinger operator is easy, it is sufficient to inspect the periodic and antiperiodic solutions. People often do the same when determining the band edges in higher dimensions through extrema of the dispersion functions at the border of the respective Brillouin zone. Quantum graphs provide a warning: there are examples of a periodic graph in which (some) band edges correspond to \emph{internal points} of the dual cell \cite{HKSW07}. Moreover, the same can happen even for graphs periodic in one dimension only \cite{EKW10}: for an array of copies of a compact graph linked by $N$ edges, as sketched in Fig.~16,
 % -------------- %
\begin{figure}[h!]
\centering
    \includegraphics[clip, trim=1cm 0cm 1cm 0cm,width=.5\textwidth]{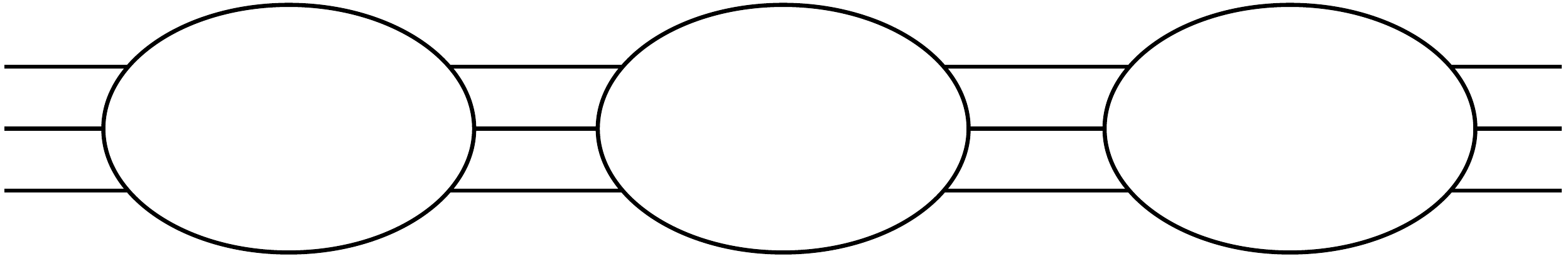}
    \caption{A graph periodic in one direction}
\end{figure}
 % -------------- %
this is can be the case for any $N\ge 2$. An example of such a graph with $N=2$ and Kirchhoff coupling worked out in \cite{EKW10} is shown in Fig.~17; note that the graph is not planar.

%\vspace{-7.5em}

 % -------------- %
\begin{figure}[h!]
\centering
    \includegraphics[clip, trim=0.5cm 0cm 0.5cm 7.5cm,angle=45,width=.6\textwidth]{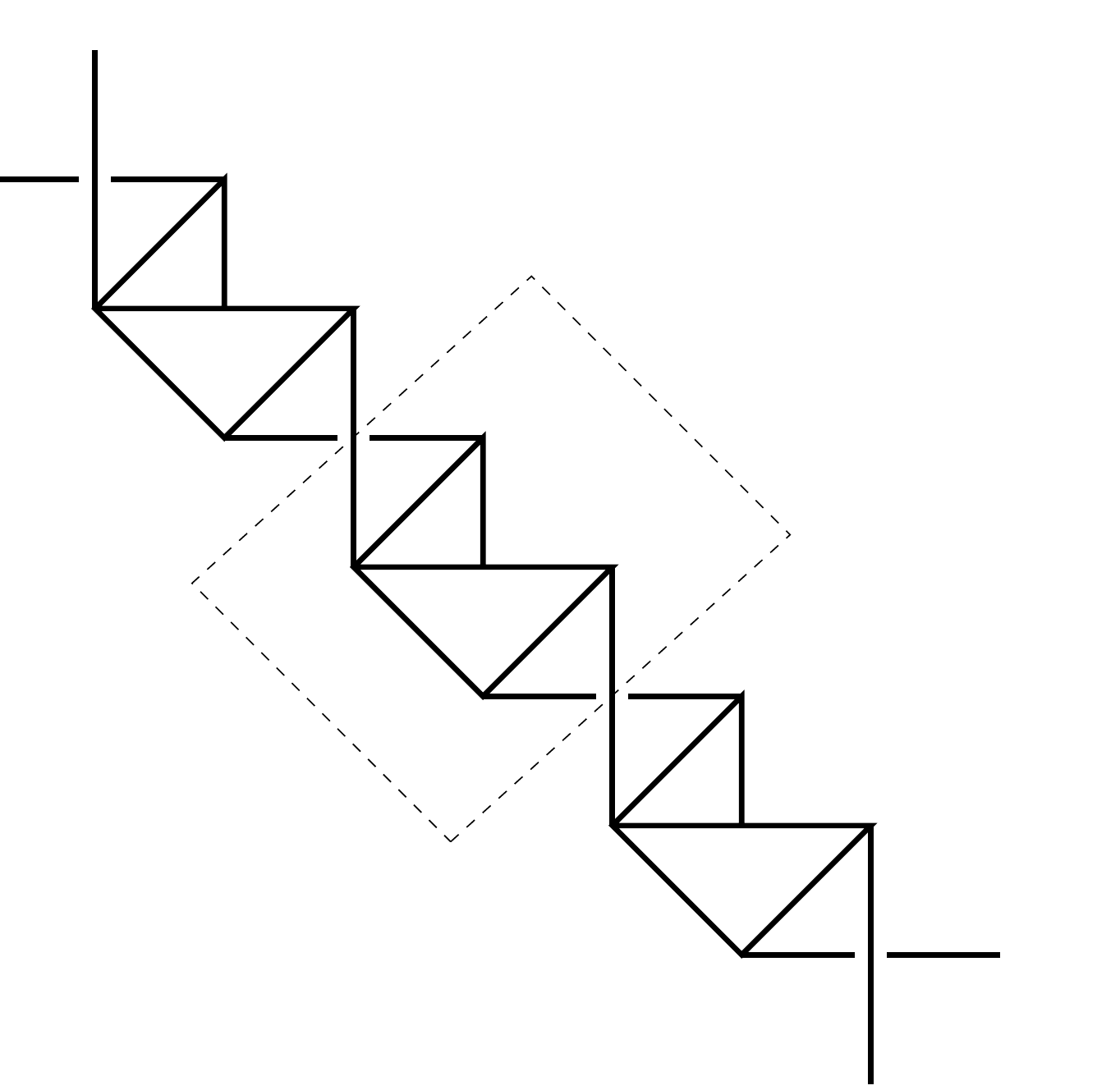}
    \vspace{-6em}
    \caption{A graph with band edges reached at points \newline inside the Brillouin zone}
\end{figure}
 % -------------- %

On the other hand, one can show that periodic and antiperiodic solutions determine the band edges if $N=1$, however, there is a catch: to demonstrate the claim one has to assume that the system is \emph{invariant with respect to time reversal}. To show that this assumption is essential, let us consider periodic \emph{comb-shaped graphs} sketched in Fig.~18 with the vertex coupling described by the matrix \eqref{rotU}.
 % -------------- %
\begin{figure}[h!]
\centering
    \includegraphics[clip, trim=1.5cm 25.5cm 7cm 2cm,width=.9\textwidth]{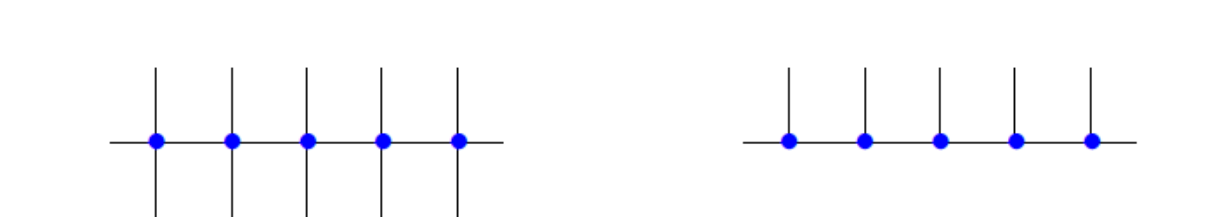}
    \caption{Comb-shaped graphs}
\end{figure}
 % -------------- %
Without going to details for which we refer to \cite{EV19} one finds that

\begin{itemize}
\setlength{\itemsep}{0pt}
\item the two-sided comb is \emph{transport-friendly}, that is, bands dominate in the high-energy regime,
\item the one-sided comb is transport-unfriendly, gaps dominate,
\item \emph{non-uniform continuity:} sending the one side edge lengths to zero in a two-sided comb does yield the one-sided comb spectrum in the sense of the set convergence in correspondence with \cite[Thm.~3.5]{BLS19}, however, in view of the above observations the convergence is highly non-uniform.
\end{itemize}

Our main question concerns the dispersion curves. The lowest ones are shown in Fig.~19; it is obvious that they reach their maxima and minima at points which are neither at the Brillouin zone boundary nor in its center.
 % -------------- %
\begin{figure}[h!]
\centering
    \includegraphics[clip, trim=0cm 0cm 0cm 0cm,width=.5\textwidth]{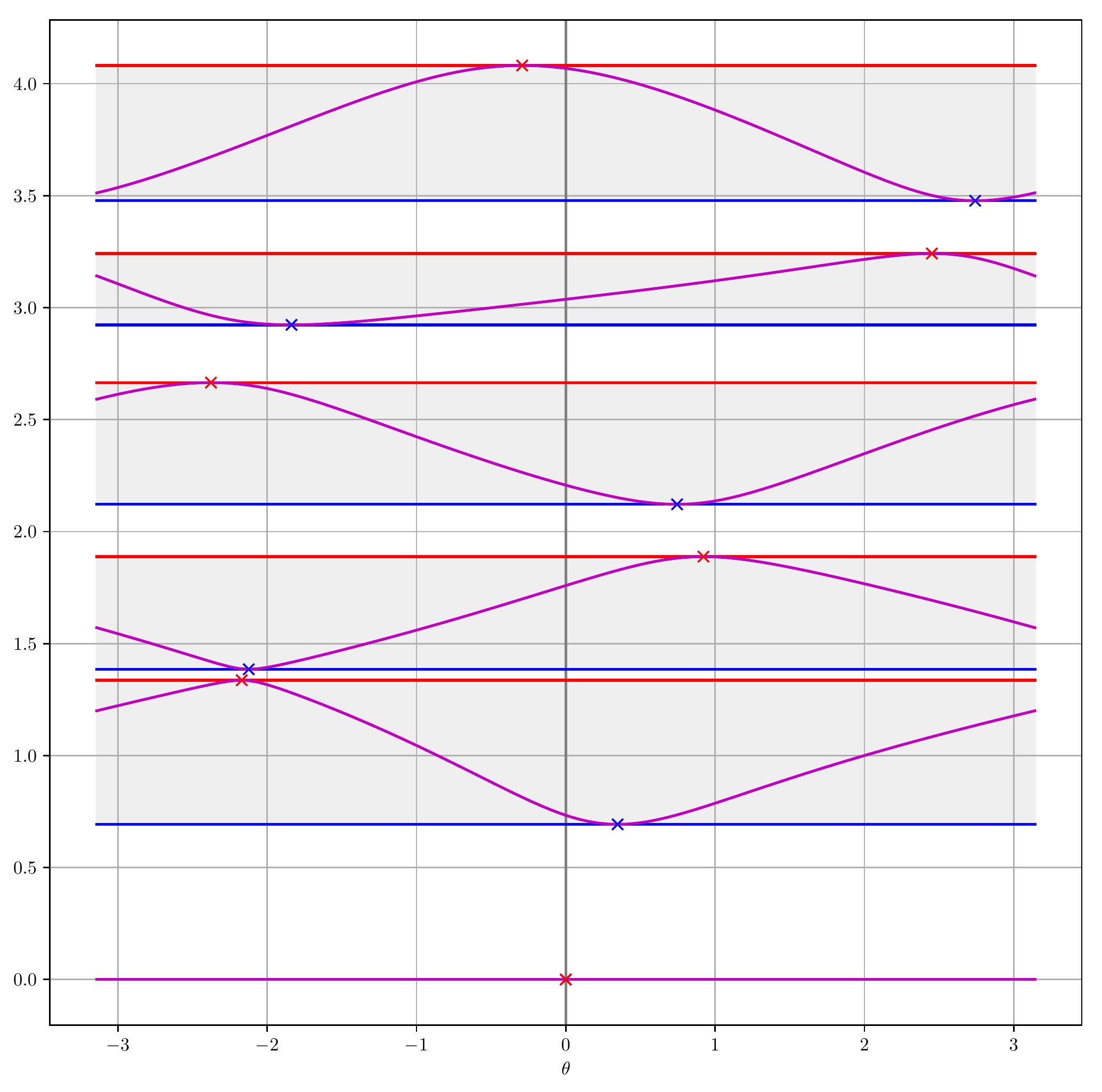}
    \caption{Dispersion curves for the comb-shaped graph}
\end{figure}
 % -------------- %

%%%%%%%%%%%%%%%%%%%%%%%%%%%%%%%%%%%%%%%%%%%%
\section{Concluding remarks} \label{s:concl}
\setcounter{equation}{0}

The above survey covered only a part of spectral problems related to quantum graphs. Other interesting questions can be addressed, for instance,

\vspace{-.7em}

\begin{itemize}

\item also for \emph{finite graphs}, the spectra of which are discrete, the topology plays often a role. An example which can be regarded as a counterpart of the results discussed in Sec.~\ref{s:timereversal} is provided by the high-energy behavior of the graphs associated with the five Platonic solids, with the vertex coupling \eqref{rotat}, where the eigenvalue asymptotics of the octahedron, as the only one with even vertex degree, differs from those of the other four \cite{EL19},

\item at the end of Section~\ref{s:intro} we mentioned the \emph{inverse problem} expressed by the modification of celebrated Mark Kac drum question \cite{vB01, GS01}. As is the PDE case, the answer is negative and one can find classes of isospectral graphs \cite{BPB09}. The topology plays again a role here, for instance through the relation between the nodal count of the graph Laplacian with Kirchhoff coupling at the vertices and the first Betti number of the corresponding graph \cite{ABB18}, see also \cite{Ba14, JH18}.

\item various problems concern the \emph{spectral optimization} with respect to graph properties. Apart from classical results such as \cite{Fr05} there is a number of recent ones \cite{BL17, BKM17, BKM19, Ro17}. As an example let us mention a generalization of the classical Ambartsumian theorem \cite{BKS18} by which the second Neumann eigenvalues on a Kirchhoff graph is minimized if the topology is trivial, that is, the edges form a segment. On the other hand, such problems become more complicated in the non-Kirchhoff case, even it is just the $\delta$ coupling \cite{EJ12},

\item for graphs with external \emph{semi-infinite edges}, often called \emph{leads}, the distribution of \emph{resonances} is of interest. It was observed by Davies and Pushnitski that the Weyl asymptotics can be violated in the Kirchhoff case if such a graph contains a \emph{balanced vertex}, in other words, the one joining the same number of external and internal edges \cite{DP11}. The result can be generalized to general couplings but it loses this simple topological formulation \cite{DEL10},

\item infinite graphs \emph{without a positive lower bound to the edge length} make the spectral analysis more involved \cite{EKMN18}, however, they allow to consider interesting cases such a graphs with a \emph{fractal} geometry \cite{AKT16},

\item other linear operators on graphs are of interest, in the first place the \emph{Dirac operator} \cite{BH03, HWK16}, also because the massless Dirac operator is a good model to describe graphene \cite{CGP09}. Other examples concerns the \emph{wave equation} \cite{FL13}, the \emph{elasticity problem} \cite {GM19, KKU15}, and one should not forget the momentum operator on graphs \cite{Ex13},

\item \emph{nonlinear operators}, in the first place nonlinear Schr\"odinger \cite{CFN17, No14} but also other nonlinear equations are of interest \cite{Ca18, Se18},

\item quantum graphs with \emph{random} parameters, both from the point of view of localization \cite{EHS07, KP08, SS00} and delocalization \cite{ASW06},

\item furthermore, one can also mention statistical properties of quantum graph spectra, in particular, the problem of \emph{quantum ergodicity} \cite{BKW04, BW16, CdV15, GKP10, ISW20, KMW03},

\item last not least, the problem of squeezing limit of Schr\"odinger operators on a \emph{Dirichlet networks} mentioned in Sec.~\ref{s:coupling} is of interest.

\end{itemize}

\vspace{-.5em}

It remains to say that the above list of open questions is by no means complete, and in addition, the references given in this paper are rather a sample than a representative collection. One nevertheless hope that this brief description could give the reader an idea how rich this field is.

%\newpage

%%%%%%%%%%%%%%%%%%%%%%%%%%%%%%
\subsection*{Acknowledgements}

The author is obliged to Prof. S.T.~Yau for the possibility to give a lecture at the 8th International Congress of Chinese Mathematicians. The research results reviewed here were supported by the Czech Science Foundation (GA\v{C}R) under Grant No.~17-01706S and by the European Union within the project CZ.02.1.01/0.0/0.0/16 019/0000778. Thanks go also to the referees whose comments helped to improve the text.

%%%%%%%%%%%%%%%%%%%%%%%%%
\bibliographystyle{99}

\begin{thebibliography}{alph}

\bibitem{ASW06}
M.~Aizenman, R.~Sims, S.~Warzel: Absolutely continuous spectra of quantum tree graphs with weak disorder, \emph{Commun. Math. Phys.} \textbf{264} (2006), 371--389.

\bibitem{ABB18}
L. Alon, R. Band, G. Berkolaiko: Nodal statistics on quantum graphs, \emph{Commun. Math. Phys.} \textbf{362} (2018), 909--948.

\bibitem{AKT16}
P.~Alonso-Ruiz, D.J.~Kelleher, A.~Teplyaev: Energy and Laplacian on Hanoi-type fractal quantum graphs, \emph{J.Phys. A: Math. Theor.} \textbf{49} (2016), 165206.

\bibitem{AJ09}
A.~Avila, S.~Jitomirskaya: The Ten Martini Problem, \emph{Ann. Math.} \textbf{170} (2009), 303--342.

\bibitem{Az64}
M.Ya.~Azbel: Energy spectrum of a conduction electron in a magnetic field, \emph{J. Exp. Theor. Phys.} \textbf{19} (1964), 634--645.

\bibitem{Ba14}
R.~Band: The nodal count $\{0,1,2,3,\dots\}$ implies the graph is a tree, \emph{Phil. Trans. Roy. Soc.} \textbf{A372} (2014), 20120504.

\bibitem{BB13}
R.~Band, G.~Berkolaiko: Universality of the momentum band density of periodic networks, \emph{Phys. Rev. Lett.} \textbf{113} (2013), 130404.

\bibitem{BL17}
R.~Band, G.~L\'evy: Quantum graphs which optimize the spectral gap, \emph{Ann. H.~Poinca\'{e}}, \textbf{18} (201), 3269--3323.

\bibitem{BPB09}
R.~Band, O.~Parzanchevski, G.~Ben-Shach: The isospectral fruits of representation theory: quantum graphs and drums, \emph{J. Phys. A: Math. Theor.} \textbf{42} (2009), 175202

\bibitem{vB01}
J.~von Below: Can one hear the shape of a network?, in \emph{Partial Differential Equations on Multistructures} (F.~Ali Mehmeti et al., eds.), M.~Dekker 2001; pp.~19-36.

\bibitem{BKW04}
G.~Berkolaiko, J.P.~Keating, B.~Winn: No quantum ergodicity for star graphs, \emph{Commun. Math. Phys.} \textbf{250} (2004), 259--285.

\bibitem{BKM17}
G.~Berkolaiko, J.B.~Kennedy, P.~Kurasov, D.~Mugnolo:
Edge connectivity and the spectral gap of combinatorial and quantum graphs, \emph{J.Phys. A: Math. Theor.} \textbf{50} (2017), 365201.

\bibitem{BKM19}
G.~Berkolaiko, J.B.~Kennedy, P.~Kurasov, D.~Mugnolo:
Surgery principles for the spectral analysis of quantum graphs, \emph{Trans. Amer. Math. Soc.} \textbf{372} (2019), 5153--5197.

\bibitem{BK13}
G. Berkolaiko, P. Kuchment: \emph{Introduction to Quantum Graphs}, AMS, Providence, R.I., 2013.

\bibitem{BLS19}
G.~Berkolaiko, Y.~Latushkin, S.~Sukhtaiev: Limits of quantum graph operators with shrinking edges, \emph{Adv. Math.} \textbf{352} (2019), 632--669.

\bibitem{BS00}
M.Sh.~Birman, T.A.~Suslina: A periodic magnetic Hamiltonian with a variable metric. The problem of absolute continuity, \emph{St. Petersburg Math. J.} \textbf{11} (2000), 203--232.

\bibitem{BH03}
J. Bolte, J.M. Harrison: Spectral statistics for the Dirac operator on graphs, \emph{J. Phys. A: Math. Gen.} \textbf{36} (2003), 2747--2769.

\bibitem{BKS18}
J.~Boman, P.~Kurasov, R.~Suhr: Schr\"odinger operators on graphs and geometry II. Spectral estimates for $L_1$-potentials and an Ambartsumian theorem, \emph{Int. Eq. Oper. Theory} \textbf{90} (2018), 40.

\bibitem{BW16}
M.~Brammall, B.~Winn: Quantum ergodicity for quantum graphs without back-scattering, \emph{Ann. H.~Poincar\'{e}} \textbf{17} (2016), 1353--1382.

\bibitem{CE07}
C.~Cacciapuoti, P.~Exner: Nontrivial edge coupling from a Dirichlet network squeezing: the case of a bent waveguide, \emph{J. Phys. A: Math. Theor.} \textbf{40} (2007), L511--L523.

\bibitem{CFN17}
C.~Cacciapuoti, D.~Finco, D.~Noja: Ground state and orbital stability for the NLS equation on a general starlike graph with potentials, \emph{Nonlinearity} \textbf{30}(2017), 3271--3303.

\bibitem{CGP09}
A.H. Castro Neto, F. Guinea, N.M.R. Peres, K.S. Novoselov, A.K. Geim: The electronic properties of graphene, \emph{Rev. Mod. Phys.} \textbf{81} (2009), 109.

\bibitem{Ca18}
M.~Cavalcante: The Korteweg-de Vries equation on a metric star graph, \emph{ZAMP} \textbf{69} (2018), 124.

\bibitem{CdV15}
Y.~Colin de Verdi\`ere: Semi-classical measures on quantum graphs and the Gauss map of the determinant manifold, \emph{Ann. H.~Poincar\'{e}} \textbf{16} (2015), 347--364.

\bibitem{DEL10}
E.B.~Davies, P.~Exner, J.~Lipovsk\'{y}: Non-Weyl asymptotics for quantum graphs with general coupling conditions, \emph{J. Phys. A: Math. Theor.} \textbf{ 43} (2010), 474013.

\bibitem{DP11}
E.B.~Davies, A.~Pushnitski: Non-Weyl resonance asymptotics for quantum graphs, \emph{Anal. PDE} \textbf{5} (2011), 729--756.

\bibitem{De+13}
C.R.~Dean et al.: Hofstadter's butterfly and the fractal quantum Hall effect in moir\'{e} superlattices, \emph{Nature} \textbf{497} (2013), 598--602.

\bibitem{Ex96a}
P.~Exner: Weakly coupled states on branching graphs, \emph{Lett. Math. Phys.} \textbf{38} (1996), 313--320.

\bibitem{Ex96b}
P.~Exner: Contact interactions on graph superlattices, \emph{J. Phys. A: Math. Gen.} \textbf{29} (1996), 87--102.

\bibitem{Ex13}
P.~Exner: Momentum operators on graphs, in ``Spectral Analysis, Differential Equations and Mathematical Physics: A Festschrift in Honor of Fritz Gesztesy's 60th Birthday'' (H.~Holden, B.~Simon, G.~Teschl, eds.), \emph{Proc. Symp. Pure Math.}, vol.~87, AMS, Providence, R.I., 2013; pp.~105--118.

\bibitem{EG96}
P.~Exner, R.~Gawlista: Band spectra of rectangular graph superlattices, \emph{Phys. Rev.} \textbf{B53} (1996), 7275--7286.

\bibitem{EHS07}
P.~Exner, M.~Helm, P.~Stollmann: Localization on a quantum graph with a random potential on the edges, \emph{Rev. Math. Phys.} \textbf{19} (2007), 923--939.

\bibitem{EJ12}
P.~Exner, M.~Jex: On the ground state of quantum graphs with attractive $\delta$-coupling, \emph{Phys. Lett.} \textbf{A376} (2012), 713--717.

\bibitem{EKMN18}
P. Exner, A. Kostenko, M. Malamud, H. Neidhardt: Spectral theory of infinite quantum graphs, \emph{Ann. H. Poincar\'{e}} \textbf{19} (2018), 3457--3510.

\bibitem{EKW10}
P.~Exner, P.~Kuchment, B.~Winn: On the location of spectral edges in $\mathbb{Z}$-periodic media, \emph{J. Phys. A: Math. Theor.} \textbf{43} (2010), 474022.

\bibitem{EL19}
P. Exner, J. Lipovsk\'{y}: Spectral asymptotics of the Laplacian on Platonic solids graphs, \emph{J. Math. Phys.} \textbf{60} (2019), 122101
\bibitem{EP05}
P.~Exner, O.~Post: Convergence of spectra of graph-like thin manifolds, \emph{J. Geom. Phys.} \textbf{54} (2005), 77--115.

\bibitem{EP09}
P. Exner, O. Post: Approximation of quantum graph vertex couplings by scaled Schroedinger operators on thin branched manifolds,
\emph{J. Phys. A: Math. Theor.} \textbf{42} (2009), 415305.

\bibitem{EP13}
P.~Exner, O.~Post: A general approximation of quantum graph vertex couplings by scaled Schr\"odinger operators on thin branched manifolds, \emph{Commun. Math. Phys.} \textbf{322} (2013), 207--227.

\bibitem{ETa18}
P.~Exner, M.~Tater: Quantum graphs with vertices of a preferred orientation, \emph{Phys. Lett.} \textbf{A382} (2018), 283--287.

\bibitem{ETu18}
P.~Exner, O.~Turek: Periodic quantum graphs from the Bethe- Sommerfeld perspective, \emph{J. Phys. A: Math. Theor.} \textbf{50} (2017), 455201.

\bibitem{ETT18}
P.~Exner, O.~Turek, M.~Tater: A family of quantum graph vertex couplings interpolating between different symmetries, \emph{J. Phys. A: Math. Theor.} \textbf{51} (2018), 285301.

\bibitem{EV17}
P.~Exner, D.~Va\v{s}ata: Cantor spectra of magnetic chain graphs, \emph{J. Phys. A: Math. Theor.} \textbf{50} (2017), 165201.

\bibitem{EV19}
P.~Exner, D.~Va\v{s}ata: Spectral properties of $\mathbb{Z}$ periodic quantum chains without time reversal invariance, \emph{in preparation}

\bibitem{FW93}
M.I.~Freidlin, A.D.~Wentzell: Diffusion processes on graphs and the averaging principle, \emph{Ann. Probab.} \textbf{21} (1993), 2215-2245.

\bibitem{FL13}
P.~Freitas, J.~Lipovsk\'{y}: Eigenvalue asymptotics for the damped wave equation on metric graphs, \emph{J. Diff. Eqs} \textbf{263} (2013), 2780--2811.

\bibitem{Fr05}
L.~Friedlander: Extremal properties of eigenvalues for a metric graph, \emph{Ann. de l'institut Fourier} \textbf{55} (2005), 199--211.

\bibitem{GKP10}
S.~Gnutzmann, J.P.~Keating, F.~Piotet:  Eigenfunction statistics on quantum graphs, \emph{Ann. Phys.} \textbf{25} (2010), 2595--2640.

\bibitem{GM19}
F.~Gregorio, D.~Mugnolo: Higher order operators on networks: hyperbolic and parabolic theory, \texttt{arXiv:1912.03297}

\bibitem{Gr08}
D.~Grieser: Spectra of graph neighborhoods and scattering, \emph{Proc. Lond. Math. Soc.} \textbf{97} (2008), 718--752.

\bibitem{GS01}
B.~Gutkin, U.~Smilansky: Can one hear the shape of a graph?, \emph{J. Phys. A: Math. Gen.} \textbf{34} (2001), 6061--6068.

\bibitem{HKSW07}
J.M.~Harrison, P.~Kuchment, A.~Sobolev, B.~Winn: On occurrence of spectral edges for periodic operators inside the Brillouin zone, \emph{J. Phys. A: Math. Theor.} \textbf{40} (2007), 7597--7618.

\bibitem{HWK16}
J.M. Harrison, T. Weyand, K. Kirsten: Zeta functions of the Dirac operator on quantum graphs, \emph{J. Math. Phys.} \textbf{57} (2016), 102301.

\bibitem{HLQZ19}
B.~Helffer, Qinghui Liu, Yanhui Qu , Qi Zhou: Positive Hausdorff dimensional spectrum for the critical almost Mathieu operator, \emph{Commun. Math. Phys.} \textbf{368} (2019), 369--382.

\bibitem{Ho76}
D.R.~Hofstadter: 1976 Energy levels and wavefunctions of Bloch electrons in rational and irrational magnetic fields \emph{Phys. Rev.} \textbf{B14}, 2239--2249.

\bibitem{ISW20}
M.~Ingremeau, M.~Sabri, B.~Winn: Quantum ergodicity for large equilateral quantum graphs, \emph{Proc. London Math. Soc.} \textbf{101} (2020), 82--109.

\bibitem{JH18}
J.S.~Juul, Ch.H.~Joyner: Isospectral discrete and quantum graphs with the same flip counts and nodal counts, \emph{J. Phys. A: Math. Theor.} \textbf{51} (2018), 245101.

\bibitem{KMW03}
J.P.~Keating, J.~Markloff, B.~Winn: Value distribution of the eigenfunctions and spectral determinants of quantum star graphs, \emph{Commun. Math. Phys.} \textbf{241} (2003), 421--452.

\bibitem{KKU15}
J.-C. Kiik, P.~Kurasov, M.~Usman: On vertex conditions for elastic systems, \emph{Phys. Lett.} \textbf{A379} (2015), 1871--1876.

\bibitem{KP08}
F.~Klopp, K.~Pankrashkin: Localization on Quantum Graphs with Random Vertex Couplings, \emph{J. Stat. Phys.} \textbf{131} (2008), 651--673.

\bibitem{KS03}
V.~Kostrykin, R.~Schrader: Quantum wires with magnetic fluxes, \emph{Commun. Math. Phys.} \textbf{237} (2003), 161--179.

\bibitem{KS97}
T.~Kottos, U.~Smilansky: Quantum chaos on graphs, \emph{Phys. Rev. Lett.} \textbf{79} (1997), 4794--4797.

\bibitem{Ku04}
P. Kuchment: Quantum graphs: II. Some spectral properties of quantum and combinatorial graphs, \emph{J. Phys. A.: Math. Gen.} \textbf{38} (2005), 4887-4900.

\bibitem{KZ01}
P. Kuchment, H. Zeng: Convergence of spectra of mesoscopic systems collapsing onto a graph, \emph{J. Math. Anal. Appl.} \textbf{258} (2001), 671--700.

\bibitem{LS16}
Y.~Last, M.~Shamis: Zero Hausdorff dimension spectrum for the almost Mathieu operator, \emph{Commun. Math. Phys.} \textbf{348} (2016), 729--750.

\bibitem{LLY12}
Y.~Lin, G\.~Lippner, S.T.~Yau: Quantum tunneling on graphs, \emph{Commun. Math. Phys.} \textbf{311} (2012), 113--132.

\bibitem{MV07}
S.~Molchanov, B.~Vainberg: Scattering solutions in networks of thin fibers: small diameter asymptotics, \emph{Commun. Math. Phys.} \textbf{273} (2007), 533--559.

\bibitem{No14}
D.~Noja: Nonlinear Schr\"odinger equation on graphs: recent results and open problems, \emph{Phil. Trans. Roy. Soc.} \textbf{A372} (2014), 20130002.

\bibitem{Pa12}
K.~Pankrashkin: Unitary dimension reduction for a class of self-adjoint extensions with applications to graph-like structures, \emph{J. Math. Anal. Appl.} \textbf{396} (2012), 640--655.

\bibitem{Pa08}
L.~Parnovski: Bethe-Sommerfeld conjecture, \emph{Ann. Henri Poincar\'{e}} \textbf{9} (2008), 457--450.

\bibitem{Pa36}
L. Pauling: The diamagnetic anisotropy of aromatic molecules, \emph{J. Chem. Phys.} \textbf{4} (1936), 673--677.

\bibitem{Po+13}
L.A.~Ponomarenko et al.: Cloning of Dirac fermions in graphene superlattices \emph{Nature} \textbf{497} (2013), 594--597.

\bibitem{Ro17}
J.~Rohleder: Eigenvalue estimates for the Laplacian on a metric tree, \emph{Proc. Amer. Math. Soc.} \textbf{145} (2017), 2119--2129.

\bibitem{RS01}
J.~Rubinstein, M.~Schatzman: Variational problems on multiply connected thin strips, {I}. {B}asic estimates and convergence of the {L}aplacian spectrum, \emph{Arch. Rat. Mech. Anal.} \textbf{160} (2001), 271-308.

\bibitem{RS53}
K. Ruedenberg, C.W. Scherr: Free-electron network model for conjugated systems. I. Theory, \emph{J. Chem. Phys.} \textbf{21} (1953), 1565--1581.

\bibitem{Sa00}
Y. Saito: The limiting equation for Neumann Laplacians on shrinking domains, \emph{El. J. Diff. Eqs}  \emph{31} (2000), 25.

\bibitem{SS00}
H.~Schanz, U.~Smilansky: Periodic-orbit theory of Anderson localization on graphs, \emph{Phys. Rev. Lett.} \textbf{84} (2000), 1427--1430.

\bibitem{SA00}
J.H.~Schenker, M.~Aizenman: The creation of spectral gaps by graph decoration, \emph{Lett. Math. Phys.} \textbf{53} (2000), 253--262.

\bibitem{Se18}
Ch.~Seifert: The linearised Korteweg-de Vries equation on general metric graphs, in ``The Diversity and Beauty of Applied Operator Theory'' (A.~B\"otcher, H.~Dym, H.~Langer, Ch.~Tretter, eds.), \emph{ Operator Theory: Advances and Applications}, vol.~268, Birkh\"auser, Basel 2018; \mbox{pp.~449--458}.

\bibitem{Sh94}
M.A.~Shubin: Discrete magnetic Laplacian, \emph{Commun. Math. Phys.} \emph{164} (1994), 259--275.

\bibitem{SK15}
P.~St\v{r}eda, J.~Ku\v{c}era: Orbital momentum and topological phase transformation, \emph{Phys. Rev.} \textbf{B92} (2015), 235152.

\bibitem{Th73}
L.E.~Thomas: Time dependent approach to scattering from impurities in a crystal,
\emph{Commun. Math.Phys.} \textbf{33} (1973), 335--343.


\end{thebibliography}

\newcommand{\etalchar}[1]{$^{#1}$}

\end{document}